\documentclass[journal]{IEEEtran}
\usepackage{amsfonts}
\usepackage{amsmath}
\usepackage{array}
\usepackage{xtab}
\usepackage{tikz}
\usepackage[europeanresistors,americaninductors,
nooldvoltagedirection]{circuitikz}
\usetikzlibrary{shapes.geometric, arrows, calc}

\allowdisplaybreaks

\usepackage{cite}
\usepackage[caption=false]{subfig}
\usepackage{xcolor}
\usepackage[normalem]{ulem}
\usepackage{soul}
\usepackage{enumitem}

\usepackage{hyperref}

\tikzset{bus/.style={fullgeneric, %
        bipoles/fullgeneric/width=0.02, bipoles/fullgeneric/height=#1
    },
    bus/.default=3
}

\newcommand{\pushright}[1]{\ifmeasuring@#1\else\omit\hfill$\displaystyle#1$\fi\ignorespaces}

\newcommand{\mycolor}{black}

\begin{document}
%
% paper title
% Titles are generally capitalized except for words such as a, an, and, as,
% at, but, by, for, in, nor, of, on, or, the, to and up, which are usually
% not capitalized unless they are the first or last word of the title.
% Linebreaks \\ can be used within to get better formatting as desired.
% Do not put math or special symbols in the title.
% \title{Invertibility Conditions for Balanced\\ Power System Admittance Matrices}
\title{Invertibility Conditions for the Admittance\\ Matrices of Balanced Power Systems}
%
%
% author names and IEEE memberships
% note positions of commas and nonbreaking spaces ( ~ ) LaTeX will not break
% a structure at a ~ so this keeps an author's name from being broken across
% two lines.
% use \thanks{} to gain access to the first footnote area
% a separate \thanks must be used for each paragraph as LaTeX2e's \thanks
% was not built to handle multiple paragraphs
%

\author{Daniel Turizo,~\IEEEmembership{Member,~IEEE}
        and Daniel K. Molzahn,~\IEEEmembership{Senior Member,~IEEE}% <-this % stops a space
% DT: I have normal membership this year at least
\thanks{Electrical and Computer Engineering, Georgia Institute of Technology, \{djturizo,molzahn\}@gatech.edu. Support from NSF contract~\#2023140.}% <-this % stops a space
}%

\maketitle

% As a general rule, do not put math, special symbols or citations
% in the abstract or keywords.
\begin{abstract}
The admittance matrix encodes the network topology and electrical parameters of a power system in order to relate the current injection and voltage phasors. Since admittance matrices are central to many power engineering analyses, their characteristics are important subjects of theoretical studies. This paper focuses on the key characteristic of \emph{invertibility}. Previous literature has presented an invertibility condition for admittance matrices. This paper first identifies and fixes a technical issue in the proof of this previously presented invertibility condition. This paper then extends this previous work by deriving new conditions that are applicable to a broader class of systems with lossless branches and transformers with off-nominal tap ratios.
\end{abstract}

% Note that keywords are not normally used for peerreview papers.

% For peer review papers, you can put extra information on the cover
% page as needed:
% \ifCLASSOPTIONpeerreview
% \begin{center} \bfseries EDICS Category: 3-BBND \end{center}
% \fi
%
% For peerreview papers, this IEEEtran command inserts a page break and
% creates the second title. It will be ignored for other modes.
\IEEEpeerreviewmaketitle

\vspace{-1em}
\section*{Notation}
%\vspace{4pt}
{\small
\begin{xtabular}{p{1cm}p{7.25cm}}
$j$   & The imaginary unit ($j^2 + 1 = 0$)  \\
$a,A$   & (No boldface letter) scalar        \\
$\bf a$ & (Boldface lowercase letter) column vector \\
$\bf A$ & (Boldface uppercase letter) matrix \\
$\mathcal{A}$ & (Calligraphic font uppercase letter) set \\
${\rm Re}(\,\cdot\,)$ & Element-wise real part operator \\
${\rm Im}(\,\cdot\,)$ & Element-wise imaginary part operator \\
$(\,\cdot\,)^*$ & Element-wise conjugate operator \\
$(\,\cdot\,)^T$ & Transpose operator \\
$(\,\cdot\,)^H$ & Conjugate transpose operator \\
${\bf 0}_{n \times m}$ & Zero matrix of size $n \times m$ \\
${\bf 0}$ & Zero matrix of appropriate size, determined from context \\
$\left\{{\bf a}\right\}_k$ & $k$-th element of vector $\bf a$ (scalar) \\
$\left\{{\bf A}\right\}_k$ & $k$-th row of matrix $\bf A$ (row vector) \\
$\left\{{\bf A}\right\}_{ik}$ & Element of matrix $\bf A$ in row $i$, column $k$ (scalar) \\
$\left|{a}\right|$ & Absolute value of scalar $a$ \\
$\left|{\mathcal{A}}\right|$ & Cardinality of set $\mathcal{A}$ \\
$\left\|{\bf a}\right\|_1$ & 1-norm of vector $\bf a$: $\left\|{\bf a}\right\|_1 = \sum\limits_k {\left| {\left\{ {\bf a} \right\}_k } \right|}$ \\
$\left\|{\bf a}\right\|$ & Euclidean norm of vector $\bf a$: $\left\|{\bf a}\right\| = (\sum\limits_k {{\left| {\left\{ {\bf a} \right\}_k } \right|}^2})^{1/2}$ \\
${\rm diag}\left({\bf a}\right)$ & Diagonal matrix such that $\left\{{{\rm diag}\left({\bf a}\right)}\right\}_{kk} = \left\{{\bf a}\right\}_k$. ${\rm diag}\left({\bf a}\right)$ has as rows and columns as the size of ${\bf a}$ \\
${\rm rank}\left({\bf A}\right)$ & Rank of matrix $\bf A$ (scalar) \\
${\rm Null}\left({\bf A}\right)$ & Null space (kernel) of matrix $\bf A$ \textcolor{\mycolor}{(The set of all vectors $\bf x$ such that ${\bf A}{\bf x} = {\bf 0}$. The null space is always a vector space.)} \\
${\rm dim}(\,\cdot\,)$ & Dimension of a vector space (scalar) \\
${\rm sym}\left({\bf B}\right)$ & Symmetric part of square matrix $\bf B$: ${\rm sym} \left({\bf B}\right) = \left({{\bf B} + {\bf B}^T}\right)/2$ \\
${\bf B} \succeq {\bf 0}$ & Square matrix $\bf B$ is positive-semidefinite (for all ${\bf x} \neq {\bf 0}$, ${\rm Re}\left( {{\bf x}^H {\bf Bx}} \right) \geq 0$), but not necessarily Hermitian \\
${\bf B} \succ {\bf 0}$ & Square matrix $\bf B$ is positive-definite (for all ${\bf x} \neq {\bf 0}$, ${\rm Re}\left( {{\bf x}^H {\bf Bx}} \right) > 0$), but not necessarily Hermitian
\end{xtabular}}

%%\pagebreak
\section{Introduction}
\IEEEPARstart{T}{he} admittance matrix, which relates the current injections to the bus voltages, is one of the most fundamental concepts in power engineering. In the phasor domain, admittance matrices are complex-valued square matrices. These matrices are used in many applications, including system modeling, power flow, optimal power flow, state estimation, stability analyses, etc. \cite{kundur, crow}.
% The admittance matrix is fundamental to apply these methods, so its properties can greatly affect the methods based on it. 
This paper thoroughly characterizes the invertibility of admittance matrices, which is a fundamental property for many power system applications.
%\IEEEPARstart{T}{he} admittance matrix, which relates the current injections to the bus voltages, is one of the more fundamental concepts in power engineering. In the phasor domain, admittance matrices are complex-valued square matrices. These matrices are used in many applications, including system modeling \cite{stevenson}, power flow \cite{crow}, economic dispatch \cite{economic_dispatch}, optimal power flow \cite{opf}, contingency analysis \cite{contingency}, state estimation \cite{state_estimation}, stability analyses \cite{sauer}, etc. This paper thoroughly characterizes the invertibility of admittance matrices, which is a fundamental property for many power system applications.

% Motivate why one should care about the admittance matrix. What is it, What is it used for, Why it is important... Circuit equations, power flow equations, stability analyses, ... Include some references -- textbooks and papers.
% One sentence giving the purpose of this paper.

Several applications directly rely on the invertibility of the admittance matrix. For instance, Kron reduction \cite{kron} is a popular technique for reducing the number of independent bus voltages modeled in a power system. The feasibility of Kron reduction is contingent on the invertibility of an appropriate sub-block of the admittance matrix. Many applications of Kron reduction assume that this procedure is feasible without performing further verification (e.g., \cite{efficient_thevenin, impedance_estimator, frequency_divider}). Additionally, various fault analysis techniques require the explicit computation of the inverse of the admittance matrix (the impedance matrix) \cite{anderson_faults}. The DC power flow \cite{stott2009} and its derivative applications \cite{ptdf, atc} also require the invertibility of admittance matrices for purely inductive systems. The invertibility of the admittance matrix is a requirement seen in both classical literature and recent research efforts (see, e.g., \cite{tenti2018, rahman2019}).
%Several applications directly rely on the invertibility of the admittance matrix. For instance, Kron reduction \cite{kron} is a popular technique for reducing the number of independent bus voltages modeled in a power system. The feasibility of Kron reduction is contingent on the invertibility of an appropriate sub-block of the admittance matrix. Additionally, various fault analysis techniques require the explicit computation of the inverse of the admittance matrix (the impedance matrix) \cite{anderson_faults}. The DC power flow \cite{dc_flow, dcbx_flow} and its derivative applications \cite{ptdf, atc} also require the invertibility of admittance matrices for purely inductive systems. The invertibility of the admittance matrix is a requirement not only seen in classical literature, but in recent research efforts as well (see, e.g., \cite{wang2018, tenti2018, rahman2019}).

% Motivate why invertibility of the admittance matrix is important. Include example applications (look at Mario's paper for examples -- Kron reduction, hybrid transmission line parameters, computing the impedance matrix for fault analyses, others?). Include references to these applications. Examples of papers where people assume invertibility of the admittance matrix -- DC power flow, DC-PTDF constructions. (Find using references to Mario's previous papers). The DC power flow 

Checking invertibility of a matrix can be accomplished via rank-revealing factorizations \cite{rrlu_fact, rrqr_fact}. However, this approach is computationally costly for large matrices. Invertibility can also be checked approximately by computing the condition number via iterative algorithms that have lower complexity than matrix factorizations \cite{cond}. However, iterative estimation of the condition number can be inaccurate \cite{cond_reliable}. In some applications, such as transmission switching \cite{transmission_switching} and topology reconfiguration \cite{reconfiguration_1, reconfiguration_2}, the admittance matrix changes as part of the problem and checking invertibility for every case is intractable. Recent research has studied the theoretical characteristics of the admittance matrix in order to guarantee invertibility without the need for computationally expensive explicit checks \cite{rank_1phase, rank_3phase, gatsis, low_theorem}. \textcolor{\mycolor}{We note that these existing theoretical results have limited applicability to practical power system models, as we discuss in Section~\ref{sec:implementation}.}
%Checking invertibility of a matrix can be accomplished via rank-revealing factorizations \cite{rrlu_fact, rrqr_fact, rrqr_sparse}. However, this approach is computationally costly for large matrices. Invertibility can also be checked approximately by computing the condition number via iterative algorithms that have lower complexity than matrix factorizations \cite{cond}. However, iterative estimation of the condition number can be inaccurate \cite{cond_reliable}. In some applications, such as transmission switching \cite{transmission_switching} and topology reconfiguration \cite{reconfiguration_1, reconfiguration_2}, the admittance matrix changes as part of the problem and checking invertibility for every case is intractable. Recent research has studied the theoretical characteristics of the admittance matrix in order to guarantee invertibility without the need for computationally expensive explicit checks \cite{rank_1phase, rank_3phase, gatsis}.

% Many papers rely on invertibility and thus rigorous conditions guaranteeing invertibility are important. Discuss Mario and Nikos Gatsis' papers and their claims. 
% Directly checking invertibility is computationally difficult (discuss computational scalability). Computational aspects of checking for invertibility, O(nbranch). Some applications where you can't directly check invertibility (transmission switching, distribution system topology reconfiguration) -- want it to be invertible for any choice of decision variables in an optimization problem where we can change topology.

One of the most important results regarding theoretical invertibility guarantees comes from~\cite{rank_1phase}. The authors of~\cite{rank_1phase} show that the admittance matrix is invertible for connected networks consisting of reciprocal branches without mutual coupling and at least one shunt element\textcolor{\mycolor}{\footnote{\textcolor{\mycolor}{A branch is said to be \emph{reciprocal} if its two-port admittance matrix is symmetrical. See~\cite{bakshi} for details.}}}. This result relies on additional modeling assumptions requiring that all admittances have positive conductances and prohibiting transformers with off-nominal tap ratios (including on-load tap changers which control the voltage magnitudes or phase shifters which control the voltage angles).

These requirements can be restrictive for practical power system models. While perfectly lossless branches do not exist in physical circuits, power system datasets often approximate certain branches as lossless. 
% Additionally, network equivalencing methods may yield datasets with negative conductances~\cite{deckmann1980studies}. 
For instance, out of the 41 systems with more than 1000 buses in the PGLib test case repository~\cite{pglib}, zero-conductance branches exist in 26 systems (63.4\%). 
% Branches with negative conductances were present in 12 systems (29.3\%).
% Imposing positive conductances in these situations to ensure admittance matrix invertibility requires carefully balancing the need for making the positive conductance values small enough to avoid significant changes to the model while large enough to ensure stability of numerical algorithms. 
We further note that transformers with off-nominal tap ratios and non-zero phase shifts are also present in many practical datasets (e.g., 39 of the aforementioned 41 PGLib systems (95.1\%)).

In addition to these modeling restrictions, there is a technical issue with the proof presented in~\cite{rank_1phase}. This paper demonstrates that the result of \cite{rank_1phase} can still be achieved and generalized to a broader class of power system models. We first detail the technical issue in the proof in~\cite{rank_1phase}. We then prove invertibility of the admittance matrix under a condition that generalizes the requirements in~\cite{rank_1phase}. The condition holds for a broad class of realistic systems, including systems with lossless branches and transformers with off-nominal tap ratios. Next we show that the theorem condition holds for networks that can be decomposed into reactive components with simple structure. Finally, we present a proof-of-concept program that implements the theorem, and we show through numerical experiments that the theorem can be applied to a wide variety of realistic power systems.

% In our proof, we generalize the results of~\cite{rank_1phase} by lifting the strict-lossiness and reciprocity requirements, so our results can be applied to systems with pure reactive elements, off-nominal tap ratio transformers and phase-shifting transformers. Finally, we discuss the implications of the new conditions and how they can be used to enforce invertibility of the admittance matrix by applying small modifications to the network model.
%This paper demonstrates that the results of \cite{rank_1phase} can still be achieved, with some small modifications. First, we detail the flaws of the proof in~\cite{rank_1phase} and provide an illustrative counterexample. Then, we prove invertibility of the admittance matrix under some additional conditions that typically hold for realistic systems. In our proof, we generalize the results of~\cite{rank_1phase} by lifting the reciprocity requirement such that our results apply to systems with phase-shifting transformers. Finally, we discuss the implications of the additional conditions and how they can be used to enforce invertibility of the admittance matrix by applying small modifications to the network model.

% Describe your contributions:
% - Identified the flaw in previous papers along with a counterexample.
% - Proving new conditions that do not suffer from this flaw and are still widely applicable to practical power systems (not overly restrictive).

The rest of the paper is organized as follows. Section~II describes the result of previous research and the technical issue in their proof. Section~III states the modifications and additional lemmas required to amend and generalize the previous result to systems with purely reactive elements and more general transformer models. Section~IV describes the implementation and numerical experiments. Section~V concludes the paper.

\section{Claims from Previous Literature \texorpdfstring{\\}{} and Limitations}
% Using the notation of \cite{rank_1phase}, we can write the admittance matrix of a power system as follows:
Borrowing the notation of \cite{rank_1phase}, the admittance matrix is (see \cite{stevenson}):
\begin{equation}
{\bf Y}_\mathcal{N} = {\bf A}_{\mathcal{L},\mathcal{N}}^T {\bf Y}_\mathcal{L} {\bf A}_{\mathcal{L},\mathcal{N}} + {\bf Y}_\mathcal{T},
\end{equation}
where \textcolor{\mycolor}{${\bf A}_{\mathcal{L},\mathcal{N}} \in \mathbb{R}^{|\mathcal{L}| \times |\mathcal{N}|}$} is the oriented incidence matrix of the network graph\textcolor{\mycolor}{\footnote{\textcolor{\mycolor}{The oriented incidence matrix relates the admittances of each branch with the nodes of that branch. The $ij$-th entry is $0$ if branch $i$ is not connected to node $j$, otherwise the entry is $\pm 1$, and the sign depends on the orientation of the branch. The orientation of the branches is arbitrary. See~\cite{stevenson} for details.}}} (excluding ground), \textcolor{\mycolor}{${\bf Y}_\mathcal{L} = {\rm diag \left( {{\bf y}_\mathcal{L}}\right)} \in \mathbb{C}^{|\mathcal{L}| \times |\mathcal{L}|}$} is the diagonal matrix with the series admittances of each branch, and \textcolor{\mycolor}{${\bf Y}_\mathcal{T} = {\rm diag \left( {{\bf y}_\mathcal{T}}\right)} \in \mathbb{C}^{|\mathcal{N}| \times |\mathcal{N}|}$} is the diagonal matrix with the total shunt admittances at each node. $\mathcal{N}$ is the set of nodes (excluding ground) and $\mathcal{L}$ is the set of branches. Reference~\cite{rank_1phase} states the following assumption and lemmas (presented here with some minor extensions as described below):

\textbf{Assumption 1.} \textit{The branches are not electromagnetically coupled and have nonzero admittance, hence ${\bf Y}_\mathcal{L}$ is full-rank.}

\textbf{Lemma 1.} \textit{The rank of the oriented incidence matrix of a connected graph with $\left| {\mathcal{N}} \right|$ nodes, ${\bf A}_{\mathcal{L},\mathcal{N}}$, is $\left| {\mathcal{N}} \right|-1$. The vector of ones ${\bf 1}$ forms a basis of the null space of ${\bf A}_{\mathcal{L},\mathcal{N}}$.}

While the second statement regarding the basis of the null space is not included in Lemma~1 as presented in~\cite{rank_1phase}, it is a well-known characteristic of oriented incidence matrices\footnote{The sum of the elements of each row of ${\bf A}_{\mathcal{N},\mathcal{L}}$ is always zero since every row has exactly one entry of 1 and one entry of -1 with the rest of the entries equal to zero; see~\cite{algebraic_graph}.} that we will use later in this paper.

\textbf{Lemma 2.} \textit{The sum of the columns of ${\bf Y}_\mathcal{N}$ equals the transpose of the sum of its rows, which also equals the vector of shunt elements ${\bf y}_\mathcal{T}$ (see \cite{arrillaga}).}

\textbf{Lemma 3.} \textit{For any matrix ${\bf M}$, ${\rm rank}\left({{\bf M}^T {\bf M}}\right)={\rm rank}\left({{\bf M}}\right)$.}

As we will discuss shortly, \emph{Lemma~3 as stated above is incorrect}. This is the technical issue in~\cite{rank_1phase} mentioned above.

\textbf{Lemma 4.} \textit{For square matrices ${\bf N}_L$ and ${\bf N}_R$ with full rank and matching size, ${\rm rank}\left({{\bf N}_L {\bf M}}\right)={\rm rank}\left({{\bf M}}\right)={\rm rank}\left({{\bf M}{\bf N}_R }\right)$. Furthermore, ${\rm Null}\left({{\bf N}_L {\bf M}}\right)={\rm Null}\left({{\bf M}}\right)$.}

While the second statement regarding the relationship between the null spaces is not included in Lemma~4 as presented in~\cite{rank_1phase}, it is a well-known result from matrix theory.\footnote{Since the only solution of ${\bf N}_L {\bf x} = {\bf 0}$ is ${\bf x} = {\bf 0}$, we make ${\bf x} = {\bf M} {\bf z}$ for some vector ${\bf z}$ and the result follows.}

One of the main results of~\cite{rank_1phase} is the following theorem:

\textbf{Theorem 1.} \textit{If the graph $\left( {\mathcal{N},\mathcal{L}} \right)$ defines a connected network and Assumption 1 holds, then:}
\begin{equation}
{\rm rank}\left( {{\bf Y}_\mathcal{N} } \right) = \left\{ {\begin{array}{*{20}l}
   {\left| \mathcal{N} \right| - 1} & {\text{if} \; \; {\bf y}_\mathcal{T} = {\bf 0}},  \\
   {\left| \mathcal{N} \right|} & \text{otherwise}.  \\
\end{array}} \right.
\end{equation}
The authors of \cite{rank_1phase} prove Theorem 1 by cases. They first assume ${\bf y}_\mathcal{T} = {\bf 0}$ and use the fact that ${\bf Y}_\mathcal{L}$ is diagonal to write it as
\begin{equation}
{\bf Y}_\mathcal{L} = {\bf B}^T {\bf B},
\end{equation}
where ${\bf B} \in \mathbb{C}^{|\mathcal{N}| \times |\mathcal{N}|}$ is full-rank. Therefore:
\begin{subequations}
\begin{align}
{\bf Y}_\mathcal{N} &= {\bf A}_{\mathcal{L},\mathcal{N}}^T {\bf B}^T {\bf B} {\bf A}_{\mathcal{L},\mathcal{N}}, \\
{\bf Y}_\mathcal{N} &= \left({{\bf B} {\bf A}_{\mathcal{L},\mathcal{N}}}\right)^T {\bf B} {\bf A}_{\mathcal{L},\mathcal{N}}, \\
{\bf Y}_\mathcal{N} &= {\bf M}^T {\bf M},
\end{align}
\end{subequations}
where ${\bf M}={\bf B} {\bf A}_{\mathcal{L},\mathcal{N}}$. According to Lemma 1, ${\bf A}_{\mathcal{L},\mathcal{N}}$ has rank $\left| {\mathcal{N}} \right| - 1$. According to Lemma 4, ${\rm rank}\left({{\bf B} {\bf A}_{\mathcal{L},\mathcal{N}}}\right)={\rm rank}\left({{\bf A}_{\mathcal{L},\mathcal{N}}}\right)$, so ${\rm rank}({\bf M})=\left| \mathcal{N} \right| - 1$. Finally, according to Lemma 3, ${\rm rank}\left({{\bf Y}_\mathcal{N}}\right)=\left| {\mathcal{N}} \right| - 1$.

There is a technical issue in the proof of Theorem~1 resulting from the fact that Lemma 3 only holds for real-valued matrices. A complex-valued counterexample is the following:
\begin{equation}
{\bf M} = \left[ {\begin{array}{*{20}c}
   1 & 0  \\
   j & 0  \\
\end{array}} \right], \qquad {\rm rank}\left( {\bf M} \right) = 1,
\end{equation}
\begin{equation}
{\bf M}^T {\bf M} = \left[ {\begin{array}{*{20}c}
   0 & 0  \\
   0 & 0  \\
\end{array}} \right], \qquad {\rm rank}\left( {{\bf M}^T {\bf M}} \right) = 0.
\end{equation}
\textcolor{\mycolor}{However, Lemma 3 holds if we use the \textit{conjugate transpose} operator $(\cdot)^H$ instead of using the transpose operator $(\cdot)^T$ (that is, we not only need to transpose the matrix, we also need to conjugate its entries as well). The corrected lemma is stated next.}

\textbf{Lemma 3 (Corrected).} \textit{For any matrix ${\bf M}$ with complex entries, ${\rm rank}\left({{\bf M}^H {\bf M}}\right)={\rm rank}\left({\bf M}\right)$. Furthermore, ${\rm Null}\left({{\bf M}^H {\bf M}}\right)={\rm Null}\left({\bf M}\right)$.}

\textbf{Proof.} Suppose a vector $\bf z$ is in the null space of $\bf M$, then:
\begin{equation}
 {\bf 0} = {\bf M}{\bf z}, \quad \Longrightarrow \quad {\bf 0} = {\bf M}^H{\bf M}{\bf z},
\end{equation}
% M}$, then:
so $\mathbf{z}$ is also in the null space of ${\bf M}^H{\bf M}$. Moreover, suppose a vector $\bf z$ is in the null space of ${\bf M}^H{\bf M}$. Then, we have
\begin{subequations}
\begin{align}
 {\bf 0} = {\bf M}^H{\bf M}{\bf z}, \quad &\Longrightarrow \quad 0 = {\bf z}^H{\bf M}^H{\bf M}{\bf z} = \left\|{{\bf M}{\bf z}}\right\|^2 \\
 &\Longrightarrow \quad {\bf 0} = {\bf M}{\bf z},
\end{align}
\end{subequations}
so $\bf z$ is also in the null space of ${\bf M}$. In conclusion, $\bf z$ is in the null space of ${\bf M}$ if and only if it is in the null space ${\bf M}^H{\bf M}$; this means that ${\rm Null}\left({{\bf M}^H {\bf M}}\right)={\rm Null}\left({\bf M}\right)$. Now we apply the rank-nullity theorem (see \cite{linear_algebra}) to complete the proof. $\hfill\square$

With the corrected version of Lemma~3 and a modeling restriction to systems where all branches are strictly lossy (have positive conductances), we can fix the proof of Theorem~1 as stated above. More specifically, the assumptions of~\cite{rank_1phase} imply Theorem~2 stated in the next section.
%The proof of Theorem~1 in~\cite{rank_1phase} is still flawed, but the theorem can be amended under some additional conditions. To understand the flaw in the proof, we consider a counterexample provided by the circuit modeling a transformer in Fig.~\ref{fig:transformer_circuit}.
%

% \noindent Let $y_t = 1 / z_t$. The transformer's turns ratio $a_t$ is an arbitrary complex number. The transformer's admittance matrix is:
% \begin{align}
% {\bf Y}_t = \left[ {\begin{array}{*{20}c}
%   {y_t} & {-a_t y_t}  \\
%   {-a_t^* y_t} & {\left| {a_t} \right|^2 y_t}  \\
% \end{array}} \right] = y_t\, {\bf a}_t\, {\bf a}_t^H,
% \end{align}
% where ${\bf a}_t^H = \left[{1, -a_t}\right]$. If $a_t$ is purely real, then $a_t^*=a_t$ and we can model the transformer with the $\pi$ circuit in Fig.~\ref{fig:transformer_circuit_pi}~\cite{kundur}.

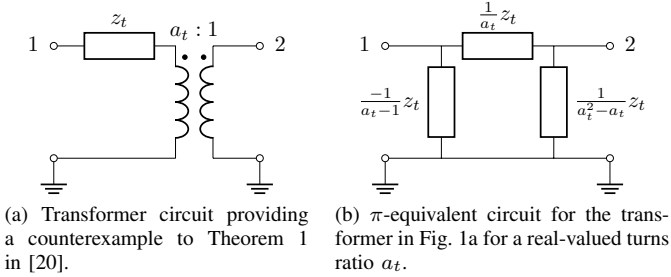
\begin{figure}[t]
\centering
\vspace{-1em}
\subfloat[Transformer circuit providing a counterexample to Theorem~1 in~\cite{rank_1phase}.]{
\scalebox{0.83}{
\begin{circuitikz}[american voltages,scale=0.6]
\draw
  % Circuit
  (-1.75,3) to [short] (-1,3)
  (-1.75,3) to [short,o-,l^=$1$] (-1.75,3)
  (-1,3) to [R,l=$z_t$] (1,3)
  (1,3) to [short] (1.5,3)
  to [L] (1.5,0) % primary winding
  node (N1) {}
  (-1.75,0) to (-1.75,0) node[ground]{}
  (1.5,0) to [short,-o] (-1.75,0)
  (2.5,0) node (N2) {}
  (2.5,0) to [L] (2.5,3)
  (3.75,0) to (3.75,0) node[ground]{}  
  (2.5,0) to [short,-o] (3.75,0)
  (2.5,3) to [short] (3,3)
  (3.75,3) to [short] (3,3)
  (3.75,3) to [short,o-,l_=$2$] (3.75,3)
  ;
    % draw transformer lines, points and tap.
  %\draw ($(N1) + (0.4,0.5)$) -- ($(N1) + (0.4,1.5)$);
  %\draw ($(N1) + (0.6,0.5)$) -- ($(N1) + (0.6,1.5)$);
  \draw [fill=black] ($(N1) + (0.25,2.7)$) circle (0.6mm);
  \draw [fill=black] ($(N1) + (0.75,2.7)$) circle (0.6mm);
  \draw ($(N1) + (0.5,3.4)$) node (tap) {$a_t:1$};
\end{circuitikz}
}
\label{fig:transformer_circuit}
}\quad
\subfloat[$\pi$-equivalent circuit for the transformer in Fig.~\ref{fig:transformer_circuit} for a real-valued turns ratio $a_t$.]{
\scalebox{0.83}{
\begin{circuitikz}[american voltages,scale=0.6]
\draw
  % Circuit
  (2,3) to [short] (3.5,3)
  (2,3) to [short,o-,l=$1$] (2,3)
  
  (3.5,3) to [R,l_=$\frac{-1}{a_t-1}z_t$] (3.5,0)
  (3.5,3) to [open,l=$\frac{1}{a_t} z_t$,yshift=0.1cm] (6.5,3)
  (3.5,3) to [R] (6.5,3)
  (6.5,3) to [R,l=$\frac{1}{a_t^2-a_t} z_t$] (6.5,0)  
  
  (8,3) to [short] (6.5,3)
  (8,3) to [short,o-,l_=$2$] (8,3)
  
  (2,0) to (2,0) node[ground]{}
  (8,0) to (8,0) node[ground]{}
  (2,0) to [short,o-o] (8,0)
  ;
\end{circuitikz}
}
\label{fig:transformer_circuit_pi}
}%
\caption{Transformer circuits.}
\vspace{-1em}
\end{figure}

We now turn our attention to the modeling restrictions of~\cite{rank_1phase}. Before generalizing Theorem~1, we need to understand why a system that violates the modeling restrictions may not satisfy the theorem. Consider the circuit modeling a transformer with an off-nominal tap ratio shown in Fig.~\ref{fig:transformer_circuit}. Let $y_t = 1 / z_t$. The transformer's turns ratio $a_t$ is an arbitrary complex number. The transformer's admittance matrix is:
\begin{align}
{\bf Y}_t = \left[ {\begin{array}{*{20}c}
  {y_t} & {-a_t y_t}  \\
  {-a_t^* y_t} & {\left| {a_t} \right|^2 y_t}  \\
\end{array}} \right] = \textcolor{\mycolor}{{\bf a}_t\, y_t\, {\bf a}_t^H}, \label{eq:single_trans}
\end{align}
where ${\bf a}_t^H = \left[{1, -a_t}\right]$. If $a_t$ is purely real, then $a_t^*=a_t$ and we can model the transformer with the $\pi$ circuit in Fig.~\ref{fig:transformer_circuit_pi}~\cite{kundur}.

The transformer's $\pi$ circuit is a two-port network with  ${\rm rank}\left( {{\bf Y}_t} \right)=1$. This $\pi$ circuit violates the requirement of strictly lossy branches if $a_t \neq 1$, as then one of the shunts will always have non-positive conductance. Notice that the impedances around the loop in the $\pi$ circuit have the sum $\frac{-1}{a_t-1}z_t + \frac{1}{a_t} z_t + \frac{1}{a_t^2-a_t} z_t = 0$. With a zero-impedance loop (i.e., a closed path through the circuit where the sum of the impedances along the path equals zero), it is mathematically possible to have non-zero voltages even in the case of zero current injections. This means that the admittance matrix is singular. More generally, admittance matrix singularity can result from other power system models with zero-impedance loops besides those associated with transformers.

The strict-lossiness restriction in~\cite{rank_1phase} requires all impedances in the power systems to have strictly positive real part. This means that the sum of the impedances over any possible loop will always have positive real part, thus being different from zero. Hence, the strict-lossiness restriction forbids the existence of zero-impedance loops. However, this also restricts the presence of transformers with off-nominal tap ratios and branches modeled as purely reactive elements, both of which appear in practical power system datasets as discussed in Section~I. To circumvent this issue, we will treat transformers as general series elements while modeling the shunt elements of the transformer $\pi$ circuit by employing an appropriate representation of the admittance matrix. \textcolor{\mycolor}{In this new representation, the branch admittances are related to the admittance matrix through a generalized version of the incidence matrix. Using this generalized incidence matrix, we can represent a transformer as a single series branch without shunts.} With this approach, the conditions we derive in this paper only forbid the existence of \textit{non-transformer} zero-impedance loops. Further, the new representation allows us to generalize Theorem~1 to systems with purely reactive elements and transformers with off-nominal tap ratios.
%Our hypothesis is that we can make Theorem 1 hold if we impose additional conditions forbidding the existence of zero-impedance loops. However, such an approach would also forbid the presence of transformers, which is unreasonably restrictive. To circumvent this issue, we will treat transformers as general series elements while modeling the shunt elements of the transformer $\pi$ circuit by employing an appropriate representation of the admittance matrix. With this approach, the conditions that we derive in this paper only forbid the existence of \textit{non-transformer} zero-impedance loops. Further, the new representation allows us to generalize Theorem~1 to systems with phase-shifting transformers.
% Further, the new representation will allow us to include phase-shifters as any other transformer, thus generalizing Theorem~1 to systems with phase-shifting transformers.

\section{Main Results}\label{sec:main_results}
This section describes our process for fixing and generalizing the main theorem. We first state and prove all necessary lemmas that will be used to prove the main result. We also declare an additional reasonable assumption that allow us to extend the result to systems with general transformer models. We then state the generalized version of the main theorem, which requires a relaxed condition in order to hold. We close this section by proving that the relaxed condition in the generalized Theorem~1 holds for power systems with reasonably common structures. \textcolor{\mycolor}{(These structures will be discussed in the conditions of Theorem 3; see Fig.~\ref{fig:alg_example} for an example system presenting these common structures).}
%This section describes our process for amending and proving the main theorem. We first state and prove all necessary lemmas that will be used to prove the main result. We also declare several additional reasonable hypotheses that allow us to generalize the result to systems with phase-shifting transformers. We then state the amended version of the main theorem, which requires some additional conditions in order to hold. We next state and prove the validity of the amended Theorem~1 under each of these conditions, individually. We end this section by discussing the implications of our results.

\subsection{Preliminaries}
We start by introducing the following assumption:

\textbf{Assumption 2.} \textit{For any series branch $l \in \mathcal{L}$ from node~$i$ to node~$k$, the admittance matrix associated with just this element can be written as ${\bf Y}_l = {\bf a}_l y_l^{\vphantom{H}} {\bf a}_l^H \in \mathbb{C}^{\left|{\mathcal{N}}\right| \times \left|{\mathcal{N}}\right|}$, where $\left\lbrace{{\bf a}_l}\right\rbrace_i = 1$, $\left\lbrace{{\bf a}_l}\right\rbrace_k = -a_l^*$ ($a_l$ is a non-zero complex number) and all other entries of ${\bf a}_l$ are zero.}

Transmission lines and transformers (including transformers with off-nominal tap ratios) satisfy Assumption 2. Transmission lines can be modeled as transformers with $a_l=1$ along with some shunt elements. This permits modeling, for instance, $\Pi$-circuit models of transmission lines. Using Assumption~2, the admittance matrix of the full system is:
\begin{equation}\label{eq:reformulated_admittance}
{\bf Y}_\mathcal{N} = \sum\limits_{l \in \mathcal{L}} {{\bf Y}_l }   + {\bf Y}_\mathcal{T}.
\end{equation}
In~\eqref{eq:reformulated_admittance}, note that ${\bf Y}_\mathcal{T}$ does not include the shunt elements in the transformers' $\pi$ circuits as these elements are instead included in ${{\bf Y}_l }$. The sum of the matrices can be rewritten as:
\begin{equation}
{\bf Y}_\mathcal{N}  = {\bf A}_{\mathcal{L},\mathcal{N}}^H {\bf Y}_\mathcal{L} {\bf A}_{\mathcal{L},\mathcal{N}} + {\bf Y}_\mathcal{T}, \label{eq:YN_standard}
\end{equation}
where, in a slight abuse of notation relative to Section~II, \textcolor{\mycolor}{${\bf A}_{\mathcal{L},\mathcal{N}} \in \mathbb{C}^{|\mathcal{L}| \times |\mathcal{N}|}$} is the \textit{generalized incidence matrix}, whose $l$-th row is $\left\lbrace{{\bf A}_{\mathcal{L},\mathcal{N}}}\right\rbrace_l = {\bf a}_l^H$; \textcolor{\mycolor}{${\bf Y}_\mathcal{L} = {\rm diag \left( {{\bf y}_\mathcal{L}}\right)} \in \mathbb{C}^{|\mathcal{L}| \times |\mathcal{L}|}$} is the diagonal matrix containing the series admittances for each branch; and \textcolor{\mycolor}{${\bf Y}_\mathcal{T} = {\rm diag \left( {{\bf y}_\mathcal{T}}\right)} \in \mathbb{C}^{|\mathcal{N}| \times |\mathcal{N}|}$} is the diagonal matrix containing the total shunt admittances at each node. \textcolor{\mycolor}{In the case of a single transformer, notice that \eqref{eq:YN_standard} reduces to \eqref{eq:single_trans} with ${\bf A}_{\mathcal{L},\mathcal{N}} = {\bf a}_t^H$, ${\bf Y}_\mathcal{L} = y_t$ and ${\bf Y}_\mathcal{T} = {\bf 0}$. In this new representation, the effect of the off-nominal tap is not represented as a shunt in ${\bf Y}_\mathcal{T}$, but is instead contained within ${\bf A}_{\mathcal{L},\mathcal{N}}$. The representation stated in \eqref{eq:YN_standard} will be the default used in the rest of the paper.}

Parallel shunts or branches with the same tap ratio can be reduced to a single branch or shunt by adding the admittances, so we assume that this reduction is always performed:

\textbf{Remark 1.} \textit{There are no parallel shunts or parallel branches with the same tap ratio.}

The connectedness condition of the network is evaluated considering its representation with parallel branches reduced. Parallel transformers with different tap ratios cannot be represented as single branch in the form stated by Assumption~2, so they are not reduced (each parallel branch individually satisfies Assumption 2, so our results are also applicable to those cases). We next state the rank-nullity theorem as we will use it several times in the paper:

\textbf{Rank-nullity theorem ([Theorem~4.4.15] in\cite{linear_algebra}).} \textit{Let ${\bf M} \in \mathbb{C}^{m \times n}$ be an arbitrary matrix, then:}
\begin{equation}
{\rm rank}\left({\bf M}\right) + {\rm dim}\left({{\rm Null}\left({\bf M}\right)}\right) = n.
\end{equation}
% \textit{see  for proof.}

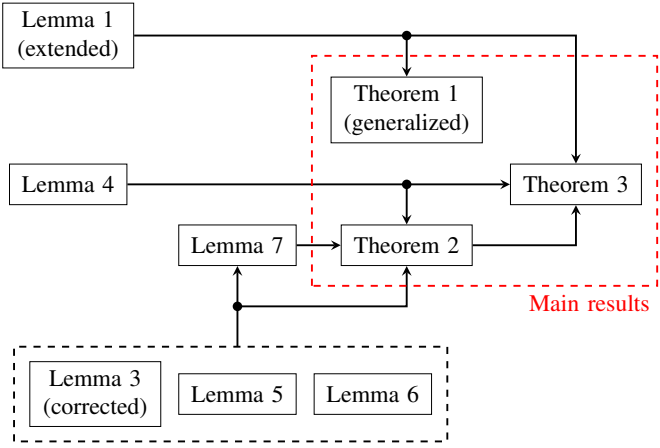
\begin{figure}[t]
\centering
\tikzstyle{process} = [rectangle, minimum width=1.0cm, minimum height=0.6cm, text centered, draw=\mycolor, text width=2.0cm]
\tikzstyle{arrow} = [thick,->,>=stealth]
\scalebox{0.9}{
\begin{tikzpicture}[color=\mycolor,x=1cm,y=1cm,node distance=1.2cm]
\node (le1) [process, text width=1.7cm] {Lemma 1 (extended)};
\node (le4) [process, text width=1.5cm, below of=le1, yshift=-1.0cm] {Lemma 4};

\node (le7) [process, text width=1.5cm, right of=le4, xshift=1.3cm, yshift=-0.9cm] {Lemma 7};

\node (le5) [process, text width=1.5cm, below of=le7, yshift=-1.0cm] {Lemma 5};
\node (le3) [process, text width=1.7cm, left of=le5, xshift=-0.9cm] {Lemma 3 (corrected)};
\node (le6) [process, text width=1.5cm, right of=le5, xshift=0.8cm] {Lemma 6};

\node (thm2) [process, text width=1.7cm, right of=le7, xshift=1.3cm] {Theorem 2};
\node (thm1) [process, above of=thm2, yshift=0.8cm] {Theorem 1 (generalized)};

\node (thm3) [process, text width=1.7cm, right of=le4, xshift=6.3cm] {Theorem 3};

\draw[\mycolor, thick, dashed] ($(le5.center)+(-3.3,-0.7)$) rectangle ($(le5.center)+(3.1,0.7)$) node[name=rec1_ne]{};

\draw [arrow,>=.] (le1) -- (le1 -| thm1) node(fork1){};
\draw [arrow] (fork1.center) -- (thm1);
\draw [arrow] (fork1.center) -| (thm3);
\draw [arrow,>=.] (le4) -- (le4 -| thm2) node(fork2){};
\draw [arrow] (fork2.center) -- (thm2);
\draw [arrow] (fork2.center) -- (thm3);
\draw [arrow] (le7) -- (thm2);
\draw [arrow] (le5.center |- rec1_ne) -- (le7);
\draw [arrow] ($0.5*(le5.center |- rec1_ne)+0.5*(le7.south)$) node(fork3){} -| (thm2);
\draw [arrow] (thm2) -| (thm3);

\fill[\mycolor] 
  (fork1) circle(2.0pt)
  (fork2) circle(2.0pt)
  (fork3) circle(2.0pt)
  ;
  
\draw[red, thick, dashed]
  ($(thm2.center)+(-1.4,-0.6)$) node[name=rec2_sw]{} rectangle ($(thm1.center -| thm3.center)+(1.2,0.8)$) node[name=rec2_ne]{}
  ($(rec2_sw -| rec2_ne)+(-1.0,0.0)$) node[anchor=north]{Main results}
  ;
\end{tikzpicture}
}
\caption{\textcolor{\mycolor}{Relationship diagram for the lemmas and theorems. An arrow going from X to Y indicates that X is used to prove Y. The three theorems are the main theoretical results of this paper.}}
\label{fig:relationships}
\vspace{-1em}
\end{figure}

\textcolor{\mycolor}{The main theoretical results of this paper are Theorems~1, 2, and 3. To prove these results, we need a series of lemmas that will be presented next. To clarify how the lemmas are related to the problem at hand, Fig.~\ref{fig:relationships} illustrates the multiple dependence relationships between the lemmas and theorems presented in this work. We start our endeavor by extending Lemma~1 to generalized incidence matrices:}

\textbf{Lemma 1 (Extended).} \textit{The rank of the generalized incidence matrix of an arbitrary connected network with $\left| {\mathcal{N}} \right|$ nodes, ${\bf A}_{\mathcal{L},\mathcal{N}} \in \mathbb{C}^{\left|{\mathcal{L}}\right| \times \left|{\mathcal{N}}\right|}$, is at least $\left| {\mathcal{N}} \right|-1$. If ${\bf A}_{\mathcal{L},\mathcal{N}}$ is not full column rank, then none of the basis vectors of its null space have null entries.}

\textbf{Proof.} Let $\mathcal{S} \subseteq \mathcal{L}$ be a set of branches forming a spanning tree of the network graph\textcolor{\mycolor}{\footnote{\textcolor{\mycolor}{To make a power engineering analogy, a \emph{spanning tree} is a subsystem obtained by removing branches from the original system until the resulting network is radial and connected. Every connected network has a spanning tree (see \cite{graph_theory}).}}}. We can order the branches of $\mathcal{L}$ by numbering all the branches of $\mathcal{S}$ first. Thus we can write ${\bf A}_{\mathcal{L},\mathcal{N}}$ in blocks as follows:
\begin{equation}
{\bf A}_{\mathcal{L},\mathcal{N}} = \left[ {\begin{array}{*{20}c}
   {\bf A}_{\mathcal{S},\mathcal{N}}  \\
   {\bf A}_{\mathcal{L} \setminus \mathcal{S},\mathcal{N}}  \\
\end{array}} \right],
\end{equation}
where ${\bf A}_{\mathcal{S},\mathcal{N}}$ is the generalized incidence matrix of the branches in $\mathcal{S}$ and ${\bf A}_{\mathcal{L} \setminus \mathcal{S},\mathcal{N}}$ is the generalized incidence matrix of the remaining branches. For any vector \textcolor{\mycolor}{$\bf x \in \mathbb{C}^{|\mathcal{N}|}$} in the null space of ${\bf A}_{\mathcal{L},\mathcal{N}}$, $\bf x$ must be orthogonal to all rows of ${\bf A}_{\mathcal{L},\mathcal{N}}$:
\begin{equation}\label{eq:nullA}
{\bf a}_s^H {\bf x} = 0, \qquad \forall s \in \mathcal{S}.
\end{equation}
Take an arbitrary branch $s$ that goes from node $i$ to node $k$, then from \eqref{eq:nullA} we have:
\begin{equation}
\left\lbrace{\bf x}\right\rbrace_i - a_s \left\lbrace{\bf x}\right\rbrace_k = 0,
\end{equation}
where $a_s$ is the tap ratio of branch $s$. We can write:
\begin{subequations}
\begin{align}
\left\lbrace{\bf x}\right\rbrace_i &= a_s \left\lbrace{\bf x}\right\rbrace_k, \\
\left\lbrace{\bf x}\right\rbrace_k &= a_s^{-1} \left\lbrace{\bf x}\right\rbrace_i.
\end{align}
\end{subequations}
We generalize this result and say that if nodes $i$ and $k$ are connected through a branch $b \in \mathcal{S}$ we can write:
\begin{equation}\label{eq:xk_from_xi}
\left\lbrace{\bf x}\right\rbrace_k = a_b^{d\left( {b,i,k} \right)} \left\lbrace{\bf x}\right\rbrace_i,
\end{equation}
where $a_b$ is the tap ratio of branch $b$, $\mathcal{S}(i,k) \subseteq \mathcal{S}$ is the (unique) set of branches in $\mathcal{S}$ forming a path from node $i$ to node $k$ (in this case the only member of $\mathcal{S}(i,k)$ is $b$), and $d\left( {b,i,k} \right)$ is a function that returns either $1$ or $-1$ depending on the direction of branch $b$ relative to the path defined by $\mathcal{S}(i,k)$ (if branch $b$ goes from node $i$ to node $k$ then $d\left( {b,i,k} \right)=-1$, otherwise $d\left( {b,i,k} \right)=1$). As $\mathcal{S}$ is a spanning tree, there exists a unique path from node 1 to every other node $k \ne 1$. Define $p_{ik}(m)$ as a function returning the node in the $m$-th position along the path from node $i$ to node $k$ ($p_{ik}(1)=i$ and $p_{ik}(1+|\mathcal{S}(i,k)|)=k$), and let $b_{ik}(m) \in \mathcal{S}_{ik}$ be the branch connecting nodes $p_{ik}(m)$ and $p_{ik}(m+1)$. Let $D_k = \left| {\mathcal{S}(1,k)} \right|$. We write $\left\lbrace{\bf x}\right\rbrace_k$ in terms of $\left\lbrace{\bf x}\right\rbrace_1$ by chaining \eqref{eq:xk_from_xi} for each pair of consecutive nodes in the path between nodes $1$ and $k$:
\begin{equation*}
1 \xrightarrow{b_{1k}(1)} p_{1k}(2) \xrightarrow{b_{1k}(2)} \,\cdots\, \xrightarrow{b_{1k}(D_k-1)} p_{1k}(D_k) \xrightarrow{b_{1k}(D_k)} k.
\end{equation*}

We backtrack the chain of equations starting from node $k$ until we reach node 1:
\begin{subequations}
\begin{align}
 \left\lbrace{\bf x}\right\rbrace_k &= a_{b_{1k}(D_k)}^{d\left( {b_{1k}(D_k), 1, k} \right)} \cdot \left\lbrace{\bf x}\right\rbrace_{p_{1k}(D_k)}, \\
\left\lbrace{\bf x}\right\rbrace_k &= a_{b_{1k}(D_k)}^{d\left( {b_{1k}(D_k), 1, k} \right)} \cdot a_{b_{1k}(D_k-1)}^{d\left( {b_{1k}(D_k-1), 1, k} \right)} \cdot \left\lbrace{\bf x}\right\rbrace_{p_{1k}(D_k-1)}, \\
 &\; \; \vdots \nonumber \\
 \left\lbrace{\bf x}\right\rbrace_k &= \left\lbrace{\bf x}\right\rbrace_1 \prod\limits_{m=1}^{D_k} {a_{b_{1k}(m)}^{d\left( {b_{1k}(m), 1, k} \right)}},
\end{align}
\end{subequations}
or written more succinctly (as the product is commutative):
\begin{equation}
\left\lbrace{\bf x}\right\rbrace_k = \left\lbrace{\bf x}\right\rbrace_1 \prod\limits_{s \in \mathcal{S}(1,k)} {a_s^{d\left( {s,1,k} \right)} }.
\end{equation}
Let $\left\lbrace{\bf x}\right\rbrace_1 = \alpha$, for an arbitrary \textcolor{\mycolor}{$\alpha \in \mathbb{C}$}. We can then write $\bf x$ as:
% Let $\left\lbrace{\bf x}\right\rbrace_1 = \alpha$ with $\alpha$ being an arbitrary value, so we can write $\bf x$ as follows:
\begin{subequations}
\begin{align}
{\bf x} &= \alpha {\bf v}, \\
\left\lbrace{\bf v}\right\rbrace_1 &= 1, \\
\left\lbrace{\bf v}\right\rbrace_k &= \prod\limits_{s \in \mathcal{S}(1,k)} {a_s^{d\left( {s,1,k} \right)} }, \qquad k = 2,\ldots,|\mathcal{N}|.
\end{align}
\end{subequations}
Since $\bf x$ has only one free parameter ($\alpha$) and ${\bf v} \neq {\bf 0}$, the rank-nullity theorem implies that ${\rm rank} \left({{\bf A}_{\mathcal{S},\mathcal{N}}}\right) = \left| {\mathcal{N}} \right|-1$. Furthermore, as $a_l \ne 0$ for all $l \in \mathcal{L}$, then all entries of $\bf v$ are non-zero.

Since $\bf x$ must also be orthogonal to all rows of ${\bf A}_{\mathcal{L} \setminus \mathcal{S},\mathcal{N}}$, we have the following equation for each row of ${\bf A}_{\mathcal{L} \setminus \mathcal{S},\mathcal{N}}$:
\begin{equation}
\alpha \left({\prod\limits_{s \in \mathcal{S}(1,i)} {a_s^{d\left( {s,1,i} \right)} } - a_l \prod\limits_{s \in \mathcal{S}(1,k)} {a_s^{d\left( {s,1,k} \right)} }}\right) = 0,
\end{equation}
for any branch $l \in \mathcal{L} \setminus \mathcal{S}$ going from node $i$ to $k$. If the term inside the parentheses is null for all rows, then the (directed) product of tap ratios $a_l$ across branches in a cycle is $1$, for all cycles. In that case, $\alpha$ is a free parameter and ${\rm rank} \left({{\bf A}_{\mathcal{L},\mathcal{N}}}\right) = \left| {\mathcal{N}} \right|-1$. Otherwise $\alpha = 0$, and so ${\rm rank} \left({{\bf A}_{\mathcal{L},\mathcal{N}}}\right) = \left| {\mathcal{N}} \right|$ (i.e., ${\bf A}_{\mathcal{L},\mathcal{N}}$ is full column rank). $\hfill\square$

We also require some new lemmas. We start with Lemma~5, which is a simple extension of Lemma 1 from \cite{gatsis}:

\textbf{Lemma 5.} \textit{Consider a matrix ${\bf Y}={\bf G}+j{\bf B} \in \mathbb{C}^{n \times n}$ with ${\bf G},{\bf B} \in \mathbb{R}^{n \times n}$. Suppose ${\bf G} \succeq {\bf 0}$, then ${\rm Null}\left({\bf Y}\right) \subseteq {\rm Null}\left({{\rm sym}\left({\bf G}\right)}\right)$ and ${\rm rank} \left({{\rm sym}\left({\bf G}\right)}\right) \le {\rm rank} \left({\bf Y}\right)$.}

\textbf{Proof.} Consider a vector ${\bf x} \in \mathbb{C}^n$ in the null space of $\bf Y$. We can write $\bf x$ in rectangular form as ${\bf x}={\bf x}_R+j {\bf x}_I$ with ${\bf x}_R,{\bf x}_I \in \mathbb{R}^n$. Using the definition of the null space, we have:
\vspace*{-1em}
\begin{subequations}
\begin{align}
 0 &= {\rm Re}\left( {{\bf x}^H {\bf Yx}} \right), \\ 
 0 &= {\bf x}_R^T {\bf Gx}_R  + {\bf x}_I^T {\bf Gx}_I  + {\bf x}_I^T {\bf Bx}_R  - {\bf x}_R^T {\bf Bx}_I.
\end{align}
\end{subequations}
The quadratic terms are real, so they only depend on the symmetric part of the matrices\footnote{For any real (possibly non-symmetric) matrix $\mathbf{A}$ and appropriately sized  real vector ${\bf x}$, the following relationships hold: ${\bf x}^T\mathbf{A}{\bf x} = ({\bf x}^T\mathbf{A}{\bf x})^T = {\bf x}^T\mathbf{A}^T{\bf x} = {\bf x}^T(\mathbf{A}/2 + \mathbf{A}^T/2){\bf x} = {\bf x}^T {\rm sym}\left( {\bf A} \right) {\bf x}$. \textcolor{\mycolor}{See \cite{quadratic_forms} for more details about symmetric quadratic forms.}}:
\begin{subequations}
\begin{align}
 0 &= {\bf x}_R^T {\rm sym}\left( {\bf G} \right){\bf x}_R  + {\bf x}_I^T {\rm sym}\left( {\bf G} \right){\bf x}_I \nonumber \\
 &\quad + {\bf x}_I^T {\rm sym}\left( {\bf B} \right){\bf x}_R  - {\bf x}_R^T {\rm sym}\left( {\bf B} \right){\bf x}_I,  \\ 
 0 &= {\bf x}_R^T {\rm sym}\left( {\bf G} \right){\bf x}_R  + {\bf x}_I^T {\rm sym}\left( {\bf G} \right){\bf x}_I \nonumber \\
 &\quad + {\bf x}_I^T {\rm sym}\left( {\bf B} \right){\bf x}_R  - {\bf x}_I^T {\rm sym}\left( {\bf B} \right){\bf x}_R, \\ 
 0 &= {\bf x}_R^T {\rm sym}\left( {\bf G} \right){\bf x}_R  + {\bf x}_I^T {\rm sym}\left( {\bf G} \right){\bf x}_I.
\end{align}
\end{subequations}
As ${\rm sym}\left( {\bf G} \right) \succeq {\bf 0}$, both terms must be non-negative. Equality only holds if both terms are zero, and hence both ${\bf x}_R$ and ${\bf x}_I$ belong to the null space of ${\rm sym}\left( {\bf G} \right)$. Therefore if ${\bf Y}{\bf x}={\bf 0}$ then ${\rm sym}\left( {\bf G} \right) {\bf x} = {\bf 0}$, so ${\rm Null}\left({\bf Y}\right) \subseteq {\rm Null}\left({{\rm sym}\left({\bf G}\right)}\right)$. We apply the rank-nullity theorem to conclude the proof. $\hfill\square$

\textbf{Lemma 6.} \textit{Let ${\bf A} \succeq {\bf 0}$ and ${\bf B} \succeq {\bf 0}$ be square matrices in $\mathbb{R}^{n \times n}$. Then the following equations hold:}
\begin{align}
&{\bf A} + {\bf B} \succeq {\bf 0}, \\
&{\rm Null}\left({{\rm sym} \left({{\bf A} + {\bf B}}\right)}\right) = \nonumber \\
&\qquad \qquad \qquad \quad \; {\rm Null}\left({{\rm sym} \left({\bf A}\right)}\right) \cap {\rm Null}\left({{\rm sym} \left({\bf B}\right)}\right), \\
&{\rm rank} \left({{\rm sym} \left({\bf A}\right)}\right), {\rm rank} \left({{\rm sym} \left({\bf B}\right)}\right) \le \nonumber \\
&\qquad \qquad \qquad \qquad \qquad \qquad \; \; \; {\rm rank} \left({{\rm sym} \left({{\bf A} + {\bf B}}\right)}\right).
\end{align}

\textbf{Proof.} Let us calculate the quadratic form of ${\bf A}+{\bf B}$:
\begin{equation}
{\bf x}^T \left({{\bf A}+{\bf B}}\right){\bf x}={\bf x}^T {\bf A}{\bf x}+{\bf x}^T {\bf B}{\bf x}.
\end{equation}
As both ${\bf A}$ and ${\bf B}$ are positive semi-definite, then ${\bf x}^T {\bf A}{\bf x}\ge 0$ and ${\bf x}^T {\bf B}{\bf x} \ge 0$, thus ${\bf x}^T \left({{\bf A}+{\bf B}}\right){\bf x} \ge 0$. Then by definition ${\bf A}+{\bf B} \succeq {\bf 0}$. Now let ${\bf z}$ be a vector in the null space of ${\rm sym} \left({{\bf A} + {\bf B}}\right)$. This implies that:
\begin{subequations}
\begin{align}
{\bf z}^T {\rm sym} \left({{\bf A} + {\bf B}}\right){\bf z}&=0, \\
{\bf z}^T {\rm sym} \left({\bf A}\right){\bf z}+{\bf z}^T {\rm sym} \left({\bf B}\right){\bf z}&=0. \label{eq:lemma6_proof}
\end{align}
\end{subequations}
Notice that ${\bf A} \succeq {\bf 0}$ implies that ${\rm sym} \left({\bf A}\right) \succeq {\bf 0}$, and similarly we have that ${\rm sym} \left({\bf B}\right) \succeq {\bf 0}$ as well. Thus both terms of \eqref{eq:lemma6_proof} are non-negative we get that
\begin{equation}
{\bf z}^T {\rm sym} \left({\bf A}\right){\bf z}=0, \qquad {\bf z}^T {\rm sym} \left({\bf B}\right){\bf z}=0,
\end{equation}
and hence ${\bf z}$ belongs to the null spaces of both ${\rm sym} \left({\bf A}\right)$ and ${\rm sym} \left({\bf B}\right)$. The converse can be proved trivially by reversing the steps, so ${\rm Null}\left({{\rm sym} \left({{\bf A} + {\bf B}}\right)}\right)={\rm Null}\left({{\rm sym} \left({\bf A}\right)}\right) \cap {\rm Null}\left({{\rm sym} \left({\bf B}\right)}\right)$. We then apply the rank-nullity theorem to conclude the proof of Lemma 6. $\hfill\square$

\textbf{Lemma 7.} \textit{Let ${\bf M} \succeq {\bf 0}$, ${\bf M} \in \mathbb{C}^{n \times n}$ be a Hermitian matrix, then ${\rm Re}\left({{\bf M}}\right) \succeq {\bf 0}$, ${\rm Null}\left({{\rm Re}\left({{\bf M}}\right)}\right) = {\rm Null}\left({\bf M}\right)$, and ${\rm rank}\left({{\rm Re}\left({{\bf M}}\right)}\right) = {\rm rank}\left({\bf M}\right)$.}

\textbf{Proof.} As $\bf M$ is Hermitian and positive-semidefinite, it can be factored as \textcolor{\mycolor}{${\bf M} = {\bf A}^H {\bf A}, {\bf A} \in \mathbb{C}^{n \times n}$}. Now we expand ${\rm Re}\left({\bf M}\right)$:
\begin{equation}
{\rm Re}\left({\bf M}\right) = {\rm Re}\left({{\bf A}}\right)^T {\rm Re}\left({{\bf A}}\right)+{\rm Im}\left({{\bf A}}\right)^T {\rm Im}\left({{\bf A}}\right).
\end{equation}
Note that each term in the right hand side is positive semidefinite, so from Lemma 6 we have that ${\rm Re}\left({\bf M}\right)$ is symmetric and positive semidefinite as well. Applying Lemmas 6 and 3:
\vspace{-1em}
\begin{subequations}
\begin{align}
{\rm Null}\left({{\rm Re}\left({\bf M}\right)}\right) &= {\rm Null}\left({{\rm Re}\left({{\bf A}}\right)}\right) \cap {\rm Null}\left({{\rm Im}\left({{\bf A}}\right)}\right), \\
{\rm Null}\left({{\rm Re}\left({\bf M}\right)}\right) &\subseteq {\rm Null}\left({{\bf A}}\right) = {\rm Null}\left({{\bf M}}\right).
\end{align}
\end{subequations}
Recall that ${\rm Re}\left({\bf M}\right)$ is symmetric and positive semidefinite, so applying Lemma~5 yields:
\begin{equation}
{\rm Null}\left({\bf M}\right) \subseteq {\rm Null}\left({{\rm Re}\left({\bf M}\right)}\right).
\end{equation}
Since ${\rm Null}\left({{\rm Re}\left({\bf M}\right)}\right) \subseteq {\rm Null}\left({{\bf M}}\right)$ and ${\rm Null}\left({\bf M}\right) \subseteq {\rm Null}\left({{\rm Re}\left({\bf M}\right)}\right)$, we have that ${\rm Null}\left({{\rm Re}\left({\bf M}\right)}\right) = {\rm Null}\left({\bf M}\right)$. The claim follows after applying the rank-nullity theorem. $\hfill\square$

\subsection{Admittance Matrix Invertibility Theorem}
We now have the tools to present the amended version of Theorem~1 and prove its validity under various conditions.

\textbf{Theorem 1 (Generalized).} \textit{Let the graph $\left( {\mathcal{N},\mathcal{L}} \right)$ define a connected network and let $\mathcal{T}$ define the shunts of the network. If Assumptions 1 and 2 hold and ${\rm Null}\left({{\bf Y}_{\mathcal{N}} }\right) \subseteq {\rm Null} \left({{\bf A}_{\mathcal{L},\mathcal{N}}}\right)$, then:}
\begin{equation}
{\rm rank}\left( {{\bf Y}_\mathcal{N} } \right) = \left\{ {\begin{array}{*{20}l}
   {{\rm rank}\left( {{\bf A}_{\mathcal{L},\mathcal{N}}} \right)} & {\text{if} \; \; \mathcal{T} = \emptyset},  \\
   {\left| \mathcal{N} \right|} & \text{otherwise}.  \\
\end{array}} \right. \label{eq:thm_rank}
\end{equation}

\textbf{Proof.}
First assume that $\mathcal{T} = \emptyset$, then
\begin{equation}
{\bf Y}_{\mathcal{N}} = {\bf A}_{\mathcal{L},\mathcal{N}}^H {\bf Y}_{\mathcal{L}} {\bf A}_{\mathcal{L},\mathcal{N}}. \label{eq:YN_noshunt}
\end{equation}
Clearly, any vector $\bf w$ such that ${\bf A}_{\mathcal{L},\mathcal{N}}\, {\bf w} = {\bf 0}$ also satisfies ${\bf Y}_{\mathcal{N}}\, {\bf w} = {\bf 0}$. This means that
\begin{equation}
{\rm Null}\left({{\bf Y}_{\mathcal{N}} }\right) \supseteq {\rm Null} \left({{\bf A}_{\mathcal{L},\mathcal{N}}}\right),
\end{equation}
so
\begin{equation}
{\rm Null}\left({{\bf Y}_{\mathcal{N}} }\right) = {\rm Null} \left({{\bf A}_{\mathcal{L},\mathcal{N}}}\right).
\end{equation}
Applying the rank-nullity theorem, we conclude that \eqref{eq:thm_rank} holds for this case. Next assume that $\mathcal{T} \neq \emptyset$, then
\begin{equation}
{\bf Y}_{\mathcal{N}} = {\bf A}_{\mathcal{L},\mathcal{N}}^H {\bf Y}_{\mathcal{L}} {\bf A}_{\mathcal{L},\mathcal{N}} + {\bf Y}_{\mathcal{T}}.
\end{equation}
If ${\bf A}_{\mathcal{L},\mathcal{N}}$ is full rank, then the fact that ${\rm Null}\left({{\bf Y}_{\mathcal{N}} }\right) \subseteq {\rm Null} \left({{\bf A}_{\mathcal{L},\mathcal{N}}}\right)$ and the rank-nullity theorem will imply that ${\bf Y}_{\mathcal{N}}$ is invertible, meaning that \eqref{eq:thm_rank} holds. If ${\bf A}_{\mathcal{L},\mathcal{N}}$ is not full rank, then we take an arbitrary vector ${\bf x} \in {\rm Null}\left({{\bf A}_{\mathcal{L},\mathcal{N}}}\right)$. From Lemma~1 (extended) we have that ${\bf x} = \alpha {\bf u}$ where $\bf u$ is a vector with no null entries. We now calculate ${\bf Y}_{\mathcal{N}}\, {\bf x}$:
\begin{subequations}
\begin{align}
{\bf Y}_{\mathcal{N}}\, {\bf x}&= \alpha {\bf Y}_{\mathcal{N}}\, {\bf u}, \\
{\bf Y}_{\mathcal{N}}\, {\bf x}&= \alpha \left({{\bf A}_{\mathcal{L},\mathcal{N}}^T {\bf Y}_{\mathcal{L}} {\bf A}_{\mathcal{L},\mathcal{N}} + {\bf Y}_{\mathcal{T}} }\right){\bf u}, \\
{\bf Y}_{\mathcal{N}}\, {\bf x}&= \alpha \left({{\bf A}_{\mathcal{L},\mathcal{N}}^T {\bf Y}_{\mathcal{L}} {\bf A}_{\mathcal{L},\mathcal{N}} {\bf u}+{\bf Y}_{\mathcal{T}} {\bf u}}\right), \\
{\bf Y}_{\mathcal{N}}\, {\bf x}&= \alpha {\bf Y}_{\mathcal{T}} {\bf u}.
\end{align}
\end{subequations}
Since ${\bf Y}_{\mathcal{T}}={\rm diag}\left({{\bf y}_{\mathcal{T}}}\right)$, ${\bf y}_{\mathcal{T}} \ne {\bf 0}$, and $\bf u$ has no null entries, we observe that ${\bf Y}_{\mathcal{N}}\, {\bf u}$ cannot be $\bf 0$ unless $ \alpha = 0$. This means the only vector in the null space of ${\bf A}_{\mathcal{L},\mathcal{N}}$ that is also in the null space of ${\bf Y}_{\mathcal{N}}$ is $\bf 0$. This implies that ${\bf Y}_{\mathcal{N}}$ is full-rank, so \eqref{eq:thm_rank} holds. $\hfill \square$

We have recovered the results of \cite{rank_1phase}, at the cost of requiring that the condition ${\rm Null}\left({{\bf Y}_{\mathcal{N}} }\right) \subseteq {\rm Null} \left({{\bf A}_{\mathcal{L},\mathcal{N}}}\right)$ holds. With the next theorem, we will show that the problem of verifying the condition for the whole network can be reduced to multiple smaller problems of the same nature.

\textbf{Theorem 2.} \textit{Let the graph $\left( {\mathcal{N},\mathcal{L}} \right)$ define a connected network and let $\mathcal{T}$ define the shunts of the network. Assumptions 1 and 2 hold, ${\rm Re}\left({y_l }\right) \ge 0$ for all $y_l$ of ${\bf y}_\mathcal{L}$ and ${\rm Re}\left({y_t}\right) \ge 0$ for all $y_t$ of ${\bf y}_\mathcal{T}$. Let $\mathcal{G}$ be the set containing the ground node, let $\mathcal{L}' \subseteq \mathcal{L}$ be the set of purely reactive branches, \textcolor{\mycolor}{let $\mathcal{N}' \subseteq \mathcal{N}$ be the set of non-isolated nodes\footnote{\textcolor{\mycolor}{An \emph{isolated node} of a graph is a node that does not have any graph branches connected to it.}} of the graph $\left({\mathcal{N}, \mathcal{L}'}\right)$}, and let $\mathcal{T}' \subseteq \mathcal{T}$ be the set of purely reactive shunts that are connected to some node in $\mathcal{N}'$. Let the reactive network $\left({\mathcal{N}', \mathcal{L}'}\right)$ have $K$ connected components ($K$ may be $0$), indexed as $\left({\mathcal{N}'(k) \cup \mathcal{G}, \mathcal{L}'(k) \cup \mathcal{T}'(k)}\right)$ for $k = 1, \ldots, K$, and let $\mathcal{T}'(k) \subseteq \mathcal{T}'$ be the set of shunts of component $k$. If the admittance matrices of all components, ${\bf Y}_{\mathcal{N}'(k)}$, satisfy that ${\rm Null}\left({{\bf Y}_{\mathcal{N}'(k)} }\right) \subseteq {\rm Null} \left({{\bf A}_{\mathcal{L}'(k),\mathcal{N}'(k)}}\right)$ (for $K=0$ this is vacuously true), then:
\begin{equation*}
{\rm Null}\left({{\bf Y}_{\mathcal{N}} }\right) \subseteq {\rm Null} \left({{\bf A}_{\mathcal{L},\mathcal{N}}}\right).
\end{equation*}
}
%\vspace{-1em}
\textcolor{\mycolor}{For the sake of clarity, we split the proof of Theorem~2 into three steps. Informally, these steps are the following:
\begin{enumerate}[label=\arabic*.]
 \item We prove that the effect of the purely reactive elements on the invertibility of the admittance matrix is independent of the elements with positive conductance. In particular, we may remove the positive conductance elements while retaining the relationship between the remaining elements and the original system.
 \item We prove that each reactive component of the system affects the invertibility of the admittance matrix independently of other components.
 \item We use steps 1 and 2 to prove the claim in Theorem~2.
\end{enumerate}
For convenience, we will also use the following conventions throughout the proof:
\begin{itemize}
  \item The incidence matrix of an empty branch set is ${\bf A}_{\emptyset,\mathcal{N}} = {\bf 0} \in \mathbb{C}^{1 \times \left|{\mathcal{N}}\right|}$.
  \item The series branch admittance matrix of an empty branch set is \textcolor{\mycolor}{${\bf Y}_{\emptyset} = 0 \in \mathbb{C}$}.
  \item The shunt admittance matrix of an empty shunt set is ${\bf Y}_{\emptyset} = {\bf 0} \in \mathbb{C}^{\left|{\mathcal{N}}\right| \times \left|{\mathcal{N}}\right|}$.
  \textcolor{\mycolor}{\item Any admittance matrix (branch, shunt, or node) has the rectangular form ${\bf Y}_{\mathcal{S}} = {\bf G}_{\mathcal{S}} + j {\bf B}_{\mathcal{S}}$. ${\bf G}_{\mathcal{S}}$ and ${\bf B}_{\mathcal{S}}$ are real matrices and $\mathcal{S}$ denotes a (branch, shunt, or node) set.}
\end{itemize}
We rely on context to discern between branch and shunt admittance matrices. 
}

\textcolor{\mycolor}{\textbf{Proof, step 1.}} Assume that both $\mathcal{L}'$ and $\mathcal{L} \setminus \mathcal{L}'$ are non-empty. We can write ${\bf A}_{\mathcal{L},\mathcal{N}}$ in block form as follows:
\begin{equation}
{\bf A}_{\mathcal{L},\mathcal{N}} = \left[ {\begin{array}{*{20}c}
   {\bf A}_{\mathcal{L}',\mathcal{N}} \\
   {\bf A}_{\mathcal{L} \setminus \mathcal{L}',\mathcal{N}} \\
\end{array}} \right]. \label{eq:AVE_blocks}
\end{equation}
Writing ${\bf Y}_{\mathcal{N}}$ in terms of the block matrices yields:
\begin{subequations}
\begin{align}
{\bf Y}_{\mathcal{N}} &= {\bf A}_{\mathcal{L}',\mathcal{N}}^H {\bf Y}_{\mathcal{L}'} {\bf A}_{\mathcal{L}',\mathcal{N}} + {\bf A}_{\mathcal{L} \setminus \mathcal{L}',\mathcal{N}}^H {\bf Y}_{\mathcal{L} \setminus \mathcal{L}'} {\bf A}_{\mathcal{L} \setminus \mathcal{L}',\mathcal{N}} \nonumber \\
&\quad + {\bf Y}_{\mathcal{T}'} + {\bf Y}_{\mathcal{T} \setminus \mathcal{T}'}, \\
{\bf Y}_{\mathcal{N}} &= {\bf A}_{\mathcal{L}',\mathcal{N}}^H {\bf Y}_{\mathcal{L}'} {\bf A}_{\mathcal{L}',\mathcal{N}} + {\bf A}_{\mathcal{L} \setminus \mathcal{L}',\mathcal{N}}^H {\bf Y}_{\mathcal{L} \setminus \mathcal{L}'} {\bf A}_{\mathcal{L} \setminus \mathcal{L}',\mathcal{N}} \nonumber \\
&\quad + {\bf Y}_{\mathcal{T}'} + {\bf G}_{\mathcal{T} \setminus \mathcal{T}'} + j {\bf B}_{\mathcal{T} \setminus \mathcal{T}'}, \label{eq:YN_L_Lp}
\end{align}
\end{subequations}
Next we compute ${\rm sym}\left( {{\bf G}_{\mathcal{N}} } \right)$ as follows:
\begin{subequations}
\begin{align}
 {\bf Y}_{\mathcal{N}} &= {\bf G}_{\mathcal{N}} + j{\bf B}_{\mathcal{N}},  \\ 
 {\bf G}_{\mathcal{N}} &= {\rm Re} \left( {{\bf A}_{\mathcal{L},\mathcal{N}}^H {\bf Y}_{\mathcal{L}} {\bf A}_{\mathcal{L},\mathcal{N}} + {\bf Y}_{\mathcal{T}}} \right), \label{eq:GN_noshunt}
\end{align}
\end{subequations}
let ${\bf A}_R$ and ${\bf A}_I$ denote the real and imaginary parts of ${\bf A}_{\mathcal{L},\mathcal{N}}$, then
\begin{subequations}
\begin{align}
{\bf G}_{\mathcal{N}} &= {\bf A}_R^T {\bf G}_{\mathcal{L}} {\bf A}_R  + {\bf A}_I^T {\bf G}_{\mathcal{L}} {\bf A}_I  + {\bf A}_I^T {\bf B}_{\mathcal{L}} {\bf A}_R  \nonumber \\
&\quad - {\bf A}_R^T {\bf B}_{\mathcal{L}} {\bf A}_I + {\bf G}_{\mathcal{T}}, \\ 
{\rm sym}\left( {{\bf G}_{\mathcal{N}} } \right) &= {\bf A}_R^T {\bf G}_{\mathcal{L}} {\bf A}_R  + {\bf A}_I^T {\bf G}_{\mathcal{L}} {\bf A}_I + {\bf G}_{\mathcal{T}}. \label{eq:Sym_GN_semipos}
\end{align}
\end{subequations}
Notice that:
\begin{subequations}
\small
\begin{align}
{\rm Re} \left({{\bf A}_{\mathcal{L},\mathcal{N}}^H {\bf G}_{\mathcal{L}} {\bf A}_{\mathcal{L},\mathcal{N}}}\right) &= {\rm Re} \left({\left({{\bf A}_R + j {\bf A}_I}\right)^H \cdot {\bf G}_{\mathcal{L}} \cdot \left({{\bf A}_R + j {\bf A}_I}\right) \phantom{{\bf A}_{\left( {{\mathcal{N}},{\mathcal{L}}} \right)}^H} \hspace*{-32pt}}\right), \\
{\rm Re} \left({{\bf A}_{\mathcal{L},\mathcal{N}}^H {\bf G}_{\mathcal{L}} {\bf A}_{\mathcal{L},\mathcal{N}}}\right) &= {\rm Re} \left({\left({{\bf A}_R^T - j {\bf A}_I^T}\right) \cdot {\bf G}_{\mathcal{L}} \cdot \left({{\bf A}_R + j {\bf A}_I}\right) \phantom{\left({{\bf A}_R^T - j {\bf A}_I^T}\right)} \hspace*{-57pt}}\right), \\
{\rm Re} \left({{\bf A}_{\mathcal{L},\mathcal{N}}^H {\bf G}_{\mathcal{L}} {\bf A}_{\mathcal{L},\mathcal{N}}}\right) &= {\bf A}_R^T {\bf G}_{\mathcal{L}} {\bf A}_R  + {\bf A}_I^T {\bf G}_{\mathcal{L}} {\bf A}_I,
\end{align}
\end{subequations}
where ${\bf A}_R^T {\bf G}_{\mathcal{L}} {\bf A}_R \succeq {\bf 0}$ and \textcolor{\mycolor}{${\bf A}_I^T {\bf G}_{\mathcal{L}} {\bf A}_I \succeq {\bf 0}$} (because all conductances are non-negative), so from Lemma 6 we have that ${\rm Re} \left({{\bf A}_{\mathcal{L},\mathcal{N}}^H {\bf G}_{\mathcal{L}} {\bf A}_{\mathcal{L},\mathcal{N}}}\right) \succeq {\bf 0}$. Replacing in \eqref{eq:Sym_GN_semipos}:
\begin{equation}\label{eq:Sym_GN_semipos_simplified}
 {\rm sym}\left( {{\bf G}_{\mathcal{N}} } \right) = {\rm Re} \left({{\bf A}_{\mathcal{L},\mathcal{N}}^H {\bf G}_{\mathcal{L}} {\bf A}_{\mathcal{L},\mathcal{N}}}\right) + {\bf G}_{\mathcal{T}}.
\end{equation}
From the definition of $\mathcal{L}'$, we have
{\setlength{\parskip}{-1em}
\begin{subequations}
%\small
\begin{align}
{\bf A}_{\mathcal{L},\mathcal{N}}^H {\bf G}_{\mathcal{L}} {\bf A}_{\mathcal{L},\mathcal{N}} &= {\bf A}_{\mathcal{L} \setminus \mathcal{L}',\mathcal{N}}^H {\bf G}_{\mathcal{L} \setminus \mathcal{L}'} {\bf A}_{\mathcal{L} \setminus \mathcal{L}',\mathcal{N}}, \\
{\rm Re} \left({{\bf A}_{\mathcal{L},\mathcal{N}}^H {\bf G}_{\mathcal{L}} {\bf A}_{\mathcal{L},\mathcal{N}}}\right) &= {\rm Re} \left({{\bf A}_{\mathcal{L} \setminus \mathcal{L}',\mathcal{N}}^H {\bf G}_{\mathcal{L} \setminus \mathcal{L}'} {\bf A}_{\mathcal{L} \setminus \mathcal{L}',\mathcal{N}}}\right), \\
{\bf G}_{\mathcal{T}} &= {\bf G}_{\mathcal{T} \setminus \mathcal{T}'}.
\end{align}
\end{subequations}
}%
Replacing in~\eqref{eq:Sym_GN_semipos_simplified} yields:
\begin{equation}
  {\rm sym}\left( {{\bf G}_{\mathcal{N}} } \right) = {\rm Re} \left({{\bf A}_{\mathcal{L} \setminus \mathcal{L}',\mathcal{N}}^H {\bf G}_{\mathcal{L} \setminus \mathcal{L}'} {\bf A}_{\mathcal{L} \setminus \mathcal{L}',\mathcal{N}}}\right) + {\bf G}_{\mathcal{T} \setminus \mathcal{T}'}. \label{eq:Sym_GN}
\end{equation}
As all conductances are non-negative, we know that ${\bf G}_{\mathcal{T} \setminus \mathcal{T}'} \succeq {\bf 0}$. Applying Lemma 6, we conclude that ${\rm sym}\left( {{\bf G}_{\mathcal{N}} } \right) \succeq {\bf 0}$ and its null space is the intersection of the null spaces of ${\rm Re} \left({{\bf A}_{\mathcal{L} \setminus \mathcal{L}',\mathcal{N}}^H {\bf G}_{\mathcal{L} \setminus \mathcal{L}'} {\bf A}_{\mathcal{L} \setminus \mathcal{L}',\mathcal{N}}}\right)$ and ${\bf G}_{\mathcal{T} \setminus \mathcal{T}'}$. As ${\rm sym}\left( {{\bf G}_{\mathcal{N}} } \right) \succeq {\bf 0}$, then ${\bf G}_{\mathcal{N}} \succeq {\bf 0}$ as well, so we can apply Lemma~5:
\begin{subequations}
\label{eq:Sym_GN_null}
\begin{align}
{\rm Null}\left({{\bf Y}_{\mathcal{N}} }\right) &\subseteq {\rm Null}\left({{\rm sym}\left( {{\bf G}_{\mathcal{N}} } \right)}\right), \\
{\rm Null}\left({{\bf Y}_{\mathcal{N}} }\right) &\subseteq {\rm Null}\left({{\rm Re} \left({{\bf A}_{\mathcal{L} \setminus \mathcal{L}',\mathcal{N}}^H {\bf G}_{\mathcal{L} \setminus \mathcal{L}'} {\bf A}_{\mathcal{L} \setminus \mathcal{L}',\mathcal{N}}}\right)}\right) \nonumber \\
&\quad \cap {\rm Null}\left({{\bf G}_{\mathcal{T} \setminus \mathcal{T}'} }\right).
\end{align}
\end{subequations}
Applying Lemma 7 yields:
\begin{align}
{\rm Null}\left({{\bf Y}_{\mathcal{N}} }\right) &\subseteq {\rm Null}\left({{\bf A}_{\mathcal{L} \setminus \mathcal{L}',\mathcal{N}}^H {\bf G}_{\mathcal{L} \setminus \mathcal{L}'} {\bf A}_{\mathcal{L} \setminus \mathcal{L}',\mathcal{N}}}\right) \nonumber \\
&\quad \cap {\rm Null}\left({{\bf G}_{\mathcal{T} \setminus \mathcal{T}'} }\right).
\end{align}
From the way $\mathcal{L}'$ is defined, we know that ${\bf G}_{\mathcal{L} \setminus \mathcal{L}'} \succ {\bf 0}$ so we can factor ${\bf G}_{\mathcal{L} \setminus \mathcal{L}'}$ as ${\bf G}_{\mathcal{L} \setminus \mathcal{L}'} = {\bf D}^H {\bf D}$, ${\bf D} \succ {\bf 0}$ (in particular, $\bf D$ is invertible). Next, we apply Lemma 3 and Lemma 4:

{\setlength{\parskip}{-1em}
\begin{subequations}
\small
\label{eq:Null_AGA}
\begin{align}
&{\rm Null}\left({{\bf A}_{\mathcal{L} \setminus \mathcal{L}',\mathcal{N}}^H {\bf G}_{\mathcal{L} \setminus \mathcal{L}'} {\bf A}_{\mathcal{L} \setminus \mathcal{L}',\mathcal{N}}}\right) =  \nonumber \\
&\hspace*{8.5em} {\rm Null}\left({\left({{\bf D} {\bf A}_{\mathcal{L} \setminus \mathcal{L}',\mathcal{N}}}\right)^H \left({{\bf D} {\bf A}_{\mathcal{L} \setminus \mathcal{L}',\mathcal{N}}}\right)}\right), \\
&{\rm Null}\left({{\bf A}_{\mathcal{L} \setminus \mathcal{L}',\mathcal{N}}^H {\bf G}_{\mathcal{L} \setminus \mathcal{L}'} {\bf A}_{\mathcal{L} \setminus \mathcal{L}',\mathcal{N}}}\right) = {\rm Null}\left({{\bf D} {\bf A}_{\mathcal{L} \setminus \mathcal{L}',\mathcal{N}}}\right), \\
&{\rm Null}\left({{\bf A}_{\mathcal{L} \setminus \mathcal{L}',\mathcal{N}}^H {\bf G}_{\mathcal{L} \setminus \mathcal{L}'} {\bf A}_{\mathcal{L} \setminus \mathcal{L}',\mathcal{N}}}\right) = {\rm Null}\left({{\bf A}_{\mathcal{L} \setminus \mathcal{L}',\mathcal{N}}}\right).
\end{align}
\end{subequations}
}%
Substituting into \eqref{eq:Sym_GN_null} yields:
\begin{equation}
{\rm Null}\left({{\bf Y}_{\mathcal{N}} }\right) \subseteq {\rm Null}\left({{\bf A}_{\mathcal{L} \setminus \mathcal{L}',\mathcal{N}}}\right) \cap {\rm Null}\left({{\bf G}_{\mathcal{T} \setminus \mathcal{T}'} }\right). \label{eq:Null_YN}
\end{equation}
With our established conventions, we note that \eqref{eq:Null_YN} holds even if $\mathcal{L} \setminus \mathcal{L}'$ is empty, so from now on we drop such assumption and only assume that $\mathcal{L}' \neq \emptyset$. Let \textcolor{\mycolor}{${\bf v} \in {\rm Null}\left({{\bf Y}_{\mathcal{N}} }\right) \subseteq \mathbb{C}^{|\mathcal{N}|}$}, then ${\bf 0} = {\bf Y}_{\mathcal{N}} {\bf v}$. From \eqref{eq:YN_L_Lp} we have that
\begin{align}
{\bf 0} &= {\bf A}_{\mathcal{L}',\mathcal{N}}^H {\bf Y}_{\mathcal{L}'} {\bf A}_{\mathcal{L}',\mathcal{N}} {\bf v} + {\bf A}_{\mathcal{L} \setminus \mathcal{L}',\mathcal{N}}^H {\bf Y}_{\mathcal{L} \setminus \mathcal{L}'} {\bf A}_{\mathcal{L} \setminus \mathcal{L}',\mathcal{N}} {\bf v} + {\bf Y}_{\mathcal{T}'} {\bf v} \nonumber \\
&\quad + {\bf G}_{\mathcal{T} \setminus \mathcal{T}'} {\bf v} + j {\bf B}_{\mathcal{T} \setminus \mathcal{T}'} {\bf v}.
\end{align}
From \eqref{eq:Null_YN}, we conclude that ${\bf G}_{\mathcal{T} \setminus \mathcal{T}'} {\bf v} = {\bf 0}$. Both ${\bf G}_{\mathcal{T} \setminus \mathcal{T}'}$ and ${\bf B}_{\mathcal{T} \setminus \mathcal{T}'}$ are diagonal, and the position of the null columns of ${\bf G}_{\mathcal{T} \setminus \mathcal{T}'}$ also correspond to null columns of ${\bf B}_{\mathcal{T} \setminus \mathcal{T}'}$. We conclude that ${\rm Null}\left({{\bf G}_{\mathcal{T} \setminus \mathcal{T}'} }\right) \subseteq {\rm Null}\left({{\bf B}_{\mathcal{T} \setminus \mathcal{T}'} }\right)$ and so ${\bf B}_{\mathcal{T} \setminus \mathcal{T}'} {\bf v} = {\bf 0}$. We also have from \eqref{eq:Null_YN} that ${\bf A}_{\mathcal{L} \setminus \mathcal{L}',\mathcal{N}} {\bf v} = {\bf 0}$. Removing these terms, the equation becomes:
\begin{equation}
{\bf 0} = {\bf A}_{\mathcal{L}',\mathcal{N}}^H {\bf Y}_{\mathcal{L}'} {\bf A}_{\mathcal{L}',\mathcal{N}} {\bf v} + {\bf Y}_{\mathcal{T}'} {\bf v}. \label{eq:YN_nullvec}
\end{equation}
\textcolor{\mycolor}{Notice that \eqref{eq:YN_nullvec} does not depend on the positive conductance elements of the system (these are the elements of $\mathcal{L} \setminus \mathcal{L}'$ and $\mathcal{T} \setminus \mathcal{T}'$). Thus, we have completed step 1.}

\textcolor{\mycolor}{\textbf{Proof, step 2.}} As we assumed that $\mathcal{L}'$ is non-empty then $K \geq 1$. We assume, without loss of generality, that the nodes of $\mathcal{N}$ and $\mathcal{N}'$ are sorted such that we can write:
\begin{subequations}
\label{eq:A_NLp}
\begin{align}
{\bf A}_{\mathcal{L}',\mathcal{N}} &= \left[ {\begin{array}{*{20}c}
   {\bf A}_{\mathcal{L}'(1),\mathcal{N}} \\
   {\vdots} \\
   {\bf A}_{\mathcal{L}'(K),\mathcal{N}} \\
\end{array}} \right], \label{eq:A_NLp_rows} \\
{\bf A}_{\mathcal{L}'(k),\mathcal{N}} &= \left[ {{\bf 0}_{|\mathcal{L}'(k)| \times |\mathcal{N}'(1)|}, \ldots, {\bf A}_{\mathcal{L}'(k),\mathcal{N}'(k)}, \ldots,}\right. \nonumber \\
&\quad\quad \left.{  {\bf 0}_{|\mathcal{L}'(k)| \times |\mathcal{N}'(K)|}, {\bf 0}_{|\mathcal{L}'(k)| \times |\mathcal{N} \setminus \mathcal{N}'|}} \right], \label{eq:A_NLpk} \\
{\bf A}_{\mathcal{L}',\mathcal{N}} &= \left[ {\begin{array}{*{20}c}
   {{\bf A}_{\mathcal{L}'(1),\mathcal{N}'(1)}} & {\cdots} & {\bf 0} & {\bf 0}\\
   {\vdots} & {\ddots} & {\vdots} & {\vdots} \\
   {\bf 0} & {\cdots} & {{\bf A}_{\mathcal{L}'(K),\mathcal{N}'(K)}} & {\bf 0} \\
\end{array}} \right].
\end{align}
\end{subequations}
Similarly:
\begin{subequations}
\label{eq:Yp_Tp_and_Y_Lp}
\begin{align}
{\bf Y}_{\mathcal{T}'} &= \left[ {\begin{array}{*{20}c}
   {{\bf Y}_{\mathcal{T}'(1)}'} & {\cdots} & {\bf 0} & {\bf 0} \\
   {\vdots} & {\ddots} & {\vdots} & {\vdots} \\
   {\bf 0} & {\cdots} & {{\bf Y}_{\mathcal{T}'(K)}'} & {\bf 0} \\
   {\bf 0} & {\cdots} & {\bf 0} & {\bf 0} \\
\end{array}} \right], \\
{\bf Y}_{\mathcal{L}'} &= \left[ {\begin{array}{*{20}c}
   {{\bf Y}_{\mathcal{L}'(1)}} & {\cdots} & {\bf 0} \\
   {\vdots} & {\ddots} & {\vdots} \\
   {\bf 0} & {\cdots} & {{\bf Y}_{\mathcal{L}'(K)}} \\
\end{array}} \right],
\end{align}
\end{subequations}
where ${{\bf Y}_{\mathcal{L}'(k)}}$ has size $|\mathcal{L}'(k)| \times |\mathcal{L}'(k)|$ and ${\bf Y}_{\mathcal{T}'(k)}'$ has size $|\mathcal{N}'(k)| \times |\mathcal{N}'(k)|$. Notice that the admittance matrix of the network $\left({\mathcal{N}'(k) \cup \mathcal{G}, \mathcal{L}'(k) \cup \mathcal{T}'(k)}\right)$, using the node in $\mathcal{G}$ as ground, is
\begin{equation}
{\bf Y}_{\mathcal{N}'(k)}  = {\bf A}_{\mathcal{L}'(k),\mathcal{N}'(k)}^H {\bf Y}_{\mathcal{L}'(k)} {\bf A}_{\mathcal{L}'(k),\mathcal{N}'(k)} + {\bf Y}_{\mathcal{T}'(k)}'. \label{eq:YNpk}
\end{equation}
Replacing \eqref{eq:A_NLp}, \eqref{eq:Yp_Tp_and_Y_Lp}, and \eqref{eq:YNpk} in \eqref{eq:YN_nullvec}, we have
\begin{equation}
{\bf 0} = \left[ {\begin{array}{*{20}c}
   {{\bf Y}_{\mathcal{N}'(1)}} & {\cdots} & {\bf 0} & {\bf 0} \\
   {\vdots} & {\ddots} & {\vdots} & {\vdots} \\
   {\bf 0} & {\cdots} & {{\bf Y}_{\mathcal{N}'(K)}} & {\bf 0} \\
   {\bf 0} & {\cdots} & {\bf 0} & {\bf 0} \\
\end{array}} \right] {\bf v}.
\end{equation}
Taking only the entries associated with nodes of $\mathcal{N}'(k)$ yields
\begin{equation}
{\bf 0} = {\bf R}_k {\bf v},
\end{equation}
where
\begin{align}
{\bf R}_k &= \left[ {{\bf 0}_{|\mathcal{N}'(k)| \times |\mathcal{N}'(1)|}, \ldots, {\bf Y}_{\mathcal{N}'(k)}, \ldots,}\right. \nonumber \\
&\quad\quad \left.{ {\bf 0}_{|\mathcal{N}'(k)| \times |\mathcal{N}'(K)|}, {\bf 0}_{|\mathcal{L}'(k)| \times |\mathcal{N} \setminus \mathcal{N}'|}} \right]. \label{eq:Rk}
\end{align}
Then, by definition:
\begin{subequations}
\begin{align}
{\bf v} &\in {\rm Null} ({\bf R}_k), \quad \forall k = 1, \ldots, K, \\
{\bf v} &\in \cap_{k=1}^{K} {\rm Null} ({\bf R}_k), \\
{\rm Null}\left({{\bf Y}_{\mathcal{N}} }\right) &\subseteq \cap_{k=1}^{K} {\rm Null} ({\bf R}_k).
\end{align}
\end{subequations}
\textcolor{\mycolor}{From \eqref{eq:Rk}, we conclude that each ${\bf R}_k$ is determined by each reactive component, independently of the others. The null space of ${\bf Y}_{\mathcal{N}}$ is contained in the null space of each ${\bf R}_k$, so we have completed step 2.}

\textcolor{\mycolor}{\textbf{Proof, step 3.}} The null space of ${\bf R}_k$ can be computed directly as the following Cartesian product:
\begin{align}
{\rm Null} ({\bf R}_k) &= \prod_{i=1}^{k-1}{\mathbb{R}^{|\mathcal{N}'(i)|}} \times {\rm Null} ({\bf Y}_{\mathcal{N}'(k)}) \times \prod_{i=k+1}^{K}{\mathbb{R}^{|\mathcal{N}'(i)|}} \nonumber \\
&\quad \times \mathbb{R}^{|\mathcal{N} \setminus \mathcal{N}'|} \label{eq:Null_Rk}
\end{align}
From the statement of Theorem 2, we have
\begin{equation}
{\rm Null} \left({{\bf Y}_{\mathcal{N}'(k)}}\right) \subseteq {\rm Null} \left({{\bf A}_{\mathcal{L}'(k),\mathcal{N}'(k)}}\right). \label{eq:Null_YNpk}
\end{equation}
Replacing this in \eqref{eq:Null_Rk} yields
\begin{align}
{\rm Null} ({\bf R}_k) &\subseteq \hspace*{-0.2em} \prod_{i=1}^{k-1}{\mathbb{R}^{|\mathcal{N}'(i)|}} \times \hspace*{-0.2em} {\rm Null} \left({{\bf A}_{\mathcal{L}'(k),\mathcal{N}'(k)}}\right) \hspace*{-0.2em} \times \hspace*{-0.7em} \prod_{i=k+1}^{K}{\mathbb{R}^{|\mathcal{N}'(i)|}} \nonumber \\
&\quad \times \mathbb{R}^{|\mathcal{N} \setminus \mathcal{N}'|}. \label{eq:Null_Rk_bound}
\end{align}
From \eqref{eq:A_NLpk}, we have
\begin{subequations}
\begin{align}
{\rm Null} \left({{\bf A}_{\mathcal{L}'(k),\mathcal{N}}}\right) &= {\rm Null} \left({\left[ {{\bf 0}_{|\mathcal{L}'(k)| \times |\mathcal{N}'(1)|}, \ldots,}\right.}\right. \nonumber \\
&\quad\quad \left.{\left.{{\bf A}_{\mathcal{L}'(k),\mathcal{N}'(k)}, \ldots, {\bf 0}_{|\mathcal{L}'(k)| \times |\mathcal{N}'(K)|}, }\right.}\right. \nonumber \\
&\quad\quad \left.{\left.{{\bf 0}_{|\mathcal{L}'(k)| \times |\mathcal{N} \setminus \mathcal{N}'|}} \right]}\right) \\
{\rm Null} \left({{\bf A}_{\mathcal{L}'(k),\mathcal{N}}}\right) &= \prod_{i=1}^{k-1}{\mathbb{R}^{|\mathcal{N}'(i)|}} \times {\rm Null} \left({{\bf A}_{\mathcal{L}'(k),\mathcal{N}'(k)}}\right) \nonumber \\
&\quad \times \prod_{i=k+1}^{K}{\mathbb{R}^{|\mathcal{N}'(i)|}} \times \mathbb{R}^{|\mathcal{N} \setminus \mathcal{N}'|}.
\end{align}
\end{subequations}
Replacing this in \eqref{eq:Null_Rk_bound} yields
\begin{equation}
{\rm Null} \left({{\bf R}_k}\right) \subseteq {\rm Null} \left({{\bf A}_{\mathcal{L}'(k),\mathcal{N}}}\right),
\end{equation}
therefore:
\begin{equation}
{\bf v} \in \cap_{k=1}^{K} {\rm Null} \left({{\bf A}_{\mathcal{L}'(k),\mathcal{N}}}\right).
\end{equation}
From \eqref{eq:A_NLp_rows}, we know that the matrices ${\bf A}_{\mathcal{L}'(k),\mathcal{N}}$ are the row blocks of ${\bf A}_{\mathcal{L}',\mathcal{N}}$, and hence
\begin{equation}
{\rm Null} \left({{\bf A}_{\mathcal{L}',\mathcal{N}}}\right) = \cap_{k=1}^{K} {\rm Null} \left({{\bf A}_{\mathcal{L}'(k),\mathcal{N}}}\right).
\end{equation}
Replacing:
\begin{subequations}
\begin{align}
{\bf v} &\in {\rm Null} \left({{\bf A}_{\mathcal{L}',\mathcal{N}}}\right), \\
{\rm Null}\left({{\bf Y}_{\mathcal{N}} }\right) &\subseteq {\rm Null} \left({{\bf A}_{\mathcal{L}',\mathcal{N}}}\right). \label{eq:Null_YN_Lp}
\end{align}
\end{subequations}
Combining \eqref{eq:Null_YN} and \eqref{eq:Null_YN_Lp} yields
\begin{subequations}
\begin{align}
{\rm Null}\left({{\bf Y}_{\mathcal{N}} }\right) &\subseteq {\rm Null} \left({{\bf A}_{\mathcal{L}',\mathcal{N}}}\right) \cap {\rm Null}\left({{\bf G}_{\mathcal{T} \setminus \mathcal{T}'} }\right) \cap {\rm Null}\left({{\bf A}_{\mathcal{L} \setminus \mathcal{L}',\mathcal{N}}}\right), \\
{\rm Null}\left({{\bf Y}_{\mathcal{N}} }\right) &\subseteq {\rm Null} \left({{\bf A}_{\mathcal{L}',\mathcal{N}}}\right) \cap {\rm Null}\left({{\bf A}_{\mathcal{L} \setminus \mathcal{L}',\mathcal{N}}}\right). \label{eq:Null_YN_ANLp_ANLmLp}
\end{align}
\end{subequations}
From \eqref{eq:AVE_blocks}, we know that ${\bf A}_{\mathcal{L}',\mathcal{N}}$ and ${\bf A}_{\mathcal{L} \setminus \mathcal{L}',\mathcal{N}}$ are the row blocks of ${\bf A}_{\mathcal{L},\mathcal{N}}$, and thus
\begin{equation}
{\rm Null} \left({{\bf A}_{\mathcal{L},\mathcal{N}}}\right) = {\rm Null} \left({{\bf A}_{\mathcal{L}',\mathcal{N}}}\right) \cap {\rm Null}\left({{\bf A}_{\mathcal{L} \setminus \mathcal{L}',\mathcal{N}}}\right). \label{eq:Null_ANL}
\end{equation}
We note that both \eqref{eq:Null_YN_ANLp_ANLmLp} and \eqref{eq:Null_ANL} hold even if $\mathcal{L}'$ is empty, so from now on we drop such assumption. Finally, from \eqref{eq:Null_YN_ANLp_ANLmLp} and \eqref{eq:Null_ANL} we conclude that in general: 
\begin{align*}
&&\hspace*{5.4em} {\rm Null}\left({{\bf Y}_{\mathcal{N}} }\right) \subseteq {\rm Null} \left({{\bf A}_{\mathcal{L},\mathcal{N}}}\right). &&\hspace*{5.4em} \square
\end{align*}

Qualitatively speaking, Theorem 2 is a recursive reduction: we can apply Theorem 1 to the network admittance matrix if we can also apply Theorem 1 to the reactive components of the network (defined by the subgraphs $\left({\mathcal{N}'(k), \mathcal{L}'(k)}\right)$). If there are no such reactive components, then we only require the standard condition of non-negative conductances in order to apply Theorem 1. We still need to prove that the conditions of Theorem 1 hold over the reactive components of the network. In the general case such proof may be too complex or even unattainable. However, we will prove that the conditions hold for common cases of reactive components with simple structures. Moreover, as we will see in the experiments section, the reactive components of practical power systems often have such structures, making the theory practically applicable. We next show the validity of Theorem 1 over several cases.

\textbf{Theorem 3.} \textit{Let the graph $\left( {\mathcal{N},\mathcal{L}} \right)$ define a connected network and let $\mathcal{T}$ define the shunts of the network. Moreover, Assumptions 1 and 2 hold. If the network satisfies at least one of the following conditions:
\begin{enumerate}[label=\arabic*)]
%  \item $\left( {\mathcal{N},\mathcal{L}} \right)$ is a tree and there exists a root node $r \in \mathcal{N}$ such that equivalent admittance of the network, measured between node $r$ and ground, is non-zero. Moreover, the equivalent admittance of any node to ground, under the condition of the parent node being grounded, is finite\textcolor{\mycolor}{\footnote{\textcolor{\mycolor}{In power system terms: a tree is a radial network, the root node is the feeder node, the parent of a node is the next node (the only one) when moving up to the feeder, a child of a node is one of the next nodes when going down from the feeder, and a leaf is one of the end nodes when going down the feeder. For formal definitions of the terms, see \cite{graph_theory}.}}}.
  \item \textcolor{\mycolor}{$\left( {\mathcal{N},\mathcal{L}} \right)$ is a tree and there exists a root node $r \in \mathcal{N}$ such that equivalent admittance of any node to ground, under the condition that the parent node (if any) is grounded, is non-zero}\textcolor{\mycolor}{\footnote{\textcolor{\mycolor}{In power system terms, a \emph{tree} is a radial network, the \emph{root node} is the feeder node, the \emph{parent} of a node is the next node (the only one) when moving up towards the feeder, a \emph{child} of a node is one of the next nodes when going down from the feeder, and a \emph{leaf} is one of the end nodes when going down the feeder. For formal definitions of the terms, see \cite{graph_theory}.}}}.
  \item $\left( {\mathcal{N},\mathcal{L}} \right)$ is a tree and $\mathcal{T} = \emptyset$.
  \item There are only inductors or there are only capacitors.
\end{enumerate}
then ${\rm Null}\left({{\bf Y}_{\mathcal{N}} }\right) \subseteq {\rm Null} \left({{\bf A}_{\mathcal{L}, \mathcal{N}}}\right)$.
}

\textbf{Proof, Condition \textit{1)}.}
Let \textcolor{\mycolor}{${\bf v} \in {\rm Null}\left({{\bf Y}_{\mathcal{N}} }\right) \subseteq \mathbb{C}^{|\mathcal{N}|}$} and define the vectors ${\bf i}_{\mathcal{L}} = {\bf Y}_{\mathcal{L}} {\bf A}_{\mathcal{L}, \mathcal{N}}\, {\bf v}$ and ${\bf i}_{\mathcal{T}} = {\bf Y}_{\mathcal{T}}\, {\bf v}$. We have that
\begin{subequations}
\begin{align}
{\bf 0} &= {\bf Y}_{\mathcal{N}}\, {\bf v} \\
{\bf 0} &= {\bf A}_{\mathcal{L}, \mathcal{N}}^H\, {\bf i}_{\mathcal{L}} + {\bf i}_{\mathcal{T}}. \label{eq:YN_loop}
\end{align}
\end{subequations}
We define $\mathcal{V}(0) \subseteq \mathcal{N}$ as the leaves of tree $\left( {\mathcal{N},\mathcal{L}} \right)$ (not including the root node $r$, see\cite{graph_theory}). For $l>0$, we define $\mathcal{V}(l) \subseteq \mathcal{N}$ as the set of nodes having all their children in $\cup_{k=0}^{l-1} \mathcal{V}(k)$, but do not belong to $\cup_{k=0}^{l-1} \mathcal{V}(k)$ themselves (i.e. $\mathcal{V}(l) \cap (\cup_{k=0}^{l-1} \mathcal{V}(k)) = \emptyset$). The \textit{height} of the tree is the unique integer $L$ such that $\{r\} = \mathcal{V}(L)$. The sets $\mathcal{V}(0), \ldots, \mathcal{V}(L)$ form a partition of $\mathcal{N}$. We also define
\begin{equation}
{\bf Y}^{i,k}  = \left[ {\begin{array}{*{20}c}
   {y_{11}^{i,k} } & {y_{12}^{i,k} }  \\
   {y_{21}^{i,k} } & {y_{22}^{i,k} }  \\
\end{array}} \right],
\end{equation}
as the $2 \times 2$ admittance matrix formed by considering only nodes $i$ and $k$ (in that order), and all branches connecting them (shunts excluded). Lastly, we define $i^{i,k}$ as
\begin{equation}
i^{i,k} = y^{i,k}_{11} \left\{{\bf v}\right\}_i + y^{i,k}_{12} \left\{{\bf v}\right\}_k.
\end{equation}
Consider a node $n \in \mathcal{V}(l)$ for some $l < L$ (so $n \neq r$). Let $\mathcal{C}(n)$ be index set of all branches connecting $n$ to some child node, let $p$ be the parent node of $n$, and let $k$ be the index of the branch connecting $n$ and $p$. The scalar equation of \eqref{eq:YN_loop} associated with node $n$ is
\begin{equation}
0 = y^{p,n}_{21} \left\{{\bf v}\right\}_p + y^{p,n}_{22} \left\{{\bf v}\right\}_n + \left\{{{\bf i}_{\mathcal{T}}}\right\}_n + \sum\nolimits_{i \in \mathcal{C}(n)} i^{n,i}. \label{eq:iL_i}
\end{equation}
We assume for induction that for any $i^{n,i} \in \mathcal{C}(n)$ we can write 
\begin{equation}
i^{n,i} = y^b_i \left\{{\bf v}\right\}_n, \qquad y^b_i \in \mathbb{C}, \label{eq:ybi}
\end{equation}
for some finite $y^b_i$. We recall that if $l = 0$ then $n$ is a leaf node, hence $\mathcal{C}(n) = \emptyset$ and the induction hypothesis holds vacuously. Let the shunt of node $n$ be $y^s_n = \left\{{{\bf Y}_{\mathcal{T}}}\right\}_{nn}$, then from \eqref{eq:ybi} and the definition of ${\bf i}_{\mathcal{T}}$ we get that
\begin{subequations}
\begin{align}
0 &= y^{p,n}_{21} \left\{{\bf v}\right\}_p + y^{p,n}_{22} \left\{{\bf v}\right\}_n + y^s_n \left\{{\bf v}\right\}_n + \sum\nolimits_{i \in \mathcal{C}(n)} y^b_i \left\{{\bf v}\right\}_n, \\
0 &= y^{p,n}_{21} \left\{{\bf v}\right\}_p + \left({y^{p,n}_{22} + y^{sb}_n}\right) \left\{{\bf v}\right\}_n, \label{eq:ysb_n_vn}
\end{align}
\end{subequations}
where
\begin{equation}
y^{sb}_n = y^s_n + \sum\nolimits_{i \in \mathcal{C}(n)} y^b_i. \label{eq:ysb_n}
\end{equation}
Multiplying by $y^{p,n}_{12}$ on both sides of \eqref{eq:ysb_n_vn} we get that
\begin{equation}
0 = y^{p,n}_{21} y^{p,n}_{12} \left\{{\bf v}\right\}_p + \left({y^{p,n}_{22} + y^{sb}_n}\right) y^{p,n}_{12} \left\{{\bf v}\right\}_n.
\end{equation}
Notice that:
\begin{subequations}
\begin{align}
i^{p,n} &= y^{p,n}_{11} \left\{{\bf v}\right\}_p + y^{p,n}_{12} \left\{{\bf v}\right\}_n, \\
y^{p,n}_{12} \left\{{\bf v}\right\}_n &= i^{p,n} - y^{p,n}_{11} \left\{{\bf v}\right\}_p,
\end{align}
\end{subequations}
hence
\begin{equation}
0 = y^{p,n}_{21} y^{p,n}_{12} \left\{{\bf v}\right\}_p + \left({y^{p,n}_{22} + y^{sb}_n}\right) \left({i^{p,n} - y^{p,n}_{11} \left\{{\bf v}\right\}_p}\right).
\end{equation}
The term $y^{p,n}_{22} + y^{sb}_n$ is the equivalent admittance between node $n$ and ground, under the condition of node $p$ being grounded. Hence $y^{p,n}_{22} + y^{sb}_n \neq 0$ according to Condition \textit{1)}, and
\begin{subequations}
\begin{align}
i^{p,n} &= \left({y^{p,n}_{11} - \frac{y^{p,n}_{12}y^{p,n}_{21}}{y^{p,n}_{22} + y^{sb}_n}}\right) \left\{{\bf v}\right\}_p, \\
i^{p,n} &= y^b_n \left\{{\bf v}\right\}_p,
\end{align}
\end{subequations}
where
\begin{equation}
y^b_n = y^{p,n}_{11} - \frac{y^{p,n}_{12}y^{p,n}_{21}}{y^{p,n}_{22} + y^{sb}_n}.
\end{equation}
We conclude that the induction hypothesis holds for any node $n \neq r$. We remark that $y^b_n$ is finite because $y^{p,n}_{22} + y^{sb}_n \neq 0$. Now we write the scalar equation of \eqref{eq:YN_loop} associated with the root node $r$:
\begin{subequations}
\begin{align}
0 &= \left\{{{\bf i}_{\mathcal{T}}}\right\}_{r} + \sum\nolimits_{i \in \mathcal{C}(r)} i^{r,i}, \\
0 &= y^s_r \left\{{\bf v}\right\}_r + \sum\nolimits_{i \in \mathcal{C}(r)} y^b_i \left\{{\bf v}\right\}_r, \\
0 &= y^{sb}_r \left\{{\bf v}\right\}_r,
\end{align}
\end{subequations}
\textcolor{\mycolor}{so $y^{sb}_r$ is the equivalent admittance of the root node $r$ to ground. As node $r$ has no parent, Condition \textit{1)} states that $y^{sb}_r$ is non-zero, so we conclude that $\left\{{\bf v}\right\}_r = 0$.} Now we propose a backward induction hypothesis: for every node $m \in \cup_{k=l}^{L} \mathcal{V}(k), l>0$ we have that $\left\{{\bf v}\right\}_m = 0$ (which trivially holds for $l=L$). We take any node $n \in \mathcal{V}(l-1)$, let $p$ be the parent node of $n$, then $p \in \cup_{k=l}^{L} \mathcal{V}(k)$ and so $\left\{{{\bf v}}\right\}_p = 0$. As $y^{p,n}_{21}$ is finite, we get from \eqref{eq:ysb_n_vn} that
\begin{equation}
0 = \left({y^{p,n}_{22} + y^{sb}_n}\right) \left\{{\bf v}\right\}_n,
\end{equation}
and as $y^{p,n}_{22} + y^{sb}_n \neq 0$ we conclude that $\left\{{\bf v}\right\}_n = 0$, proving the induction hypothesis. This means that ${\bf v} = {\bf 0}$, so
\begin{align*}
&&\hspace*{3.9em} {\rm Null}\left({{\bf Y}_{\mathcal{N}} }\right) = \{{\bf 0}\} \subseteq {\rm Null} \left({{\bf A}_{\mathcal{L}, \mathcal{N}}}\right). &&\hspace*{3.9em} \square
\end{align*}

\textbf{Proof, Condition \textit{2)}.}
In this case we have that ${\bf Y}_{\mathcal{N}} = {\bf A}_{\mathcal{L}, \mathcal{N}}^H {\bf Y}_{\mathcal{L}} {\bf A}_{\mathcal{L}, \mathcal{N}}$ (see \eqref{eq:YN_noshunt}), and thus ${\rm Null}\left({{\bf Y}_{\mathcal{N}} }\right) \supseteq {\rm Null} \left({{\bf A}_{\mathcal{L}, \mathcal{N}}}\right)$. We also know, since $\left( {\mathcal{N},\mathcal{L}} \right)$ is a tree, that the network has exactly $|\mathcal{N}|-1$ branches. This means that ${\bf A}_{\mathcal{L}, \mathcal{N}}$ has size $|\mathcal{N}|-1 \times |\mathcal{N}|$, so from Lemma~1 we have that ${\rm rank} \left({{\bf A}_{\mathcal{L}, \mathcal{N}}}\right) = |\mathcal{N}| - 1$. Applying the Frobenius inequality (see exercise 4.5.17 in \cite{linear_algebra}) to \eqref{eq:YN_noshunt}, we have
\begin{subequations}
\begin{align}
&{\rm rank} \left({{\bf A}_{\mathcal{L}, \mathcal{N}}^H {\bf Y}_{\mathcal{L}}}\right) + {\rm rank} \left({{\bf Y}_{\mathcal{L}} {\bf A}_{\mathcal{L}, \mathcal{N}}}\right) \leq \nonumber \\
&\hspace*{6.5em} {\rm rank} \left({{\bf Y}_{\mathcal{L}}}\right) + {\rm rank} \left({{\bf A}_{\mathcal{L}, \mathcal{N}}^H {\bf Y}_{\mathcal{L}} {\bf A}_{\mathcal{L}, \mathcal{N}}}\right), \\
&{\rm rank} \left({{\bf Y}_{\mathcal{N}}}\right) \geq {\rm rank} \left({{\bf A}_{\mathcal{L}, \mathcal{N}}^H {\bf Y}_{\mathcal{L}}}\right) + {\rm rank} \left({{\bf Y}_{\mathcal{L}} {\bf A}_{\mathcal{L}, \mathcal{N}}}\right) \nonumber \\
&\hspace*{6.5em} - {\rm rank} \left({{\bf Y}_{\mathcal{L}}}\right).
\end{align}
\end{subequations}
Applying Lemma 4 and the fact that ${\bf Y}_{\mathcal{L}}$ is square and full rank, we get that
\begin{subequations}
\begin{align}
{\rm rank} \left({{\bf Y}_{\mathcal{N}}}\right) &\geq {\rm rank} \left({{\bf A}_{\mathcal{L}, \mathcal{N}}^H}\right) + {\rm rank} \left({{\bf A}_{\mathcal{L}, \mathcal{N}}}\right) \nonumber \\
&\quad - {\rm rank} \left({{\bf Y}_{\mathcal{L}}}\right), \\
{\rm rank} \left({{\bf Y}_{\mathcal{N}}}\right) &\geq |\mathcal{N}| - 1.
\end{align}
\end{subequations}
Applying the rank-nullity theorem:
\begin{equation}
{\rm dim} \left({{\rm Null}\left({{\bf Y}_{\mathcal{N}} }\right)}\right) \leq 1 = {\rm dim} \left({{\rm Null}\left({{\bf A}_{\mathcal{L}, \mathcal{N}}}\right)}\right),
\end{equation}
but ${\rm Null}\left({{\bf Y}_{\mathcal{N}} }\right) \supseteq {\rm Null} \left({{\bf A}_{\mathcal{L}, \mathcal{N}}}\right)$, and as they have equal dimension then ${\rm Null}\left({{\bf Y}_{\mathcal{N}} }\right) = {\rm Null} \left({{\bf A}_{\mathcal{L}, \mathcal{N}}}\right)$. This trivially implies that ${\rm Null}\left({{\bf Y}_{\mathcal{N}} }\right) \subseteq {\rm Null} \left({{\bf A}_{\mathcal{L}, \mathcal{N}}}\right)$. $\hfill \square$

\textbf{Proof, Condition \textit{3)}.}
As the network is purely inductive (or purely capacitive) we can write
\begin{subequations}
\begin{align}
{\bf Y}_{\mathcal{N}} &= {\bf A}_{\mathcal{L}, \mathcal{N}}^H \left({jk {\bf B}_{\mathcal{L}}}\right) {\bf A}_{\mathcal{L}, \mathcal{N}} + \left({jk {\bf B}_{\mathcal{T}}}\right), \\
{\bf Y}_{\mathcal{N}} &= jk \left({{\bf A}_{\mathcal{L}, \mathcal{N}}^H {\bf B}_{\mathcal{L}} {\bf A}_{\mathcal{L}, \mathcal{N}} + {\bf B}_{\mathcal{T}}}\right),
\end{align}
\end{subequations}
where ${\bf B}_{\mathcal{L}}$ and ${\bf B}_{\mathcal{T}}$ are diagonal matrices with non-negative real entries, and $k=1$ if the network is purely capacitive or $k=-1$ if the network is purely inductive. Now we consider an alternative network with a set of nodes $\mathcal{N}'$ identical to $\mathcal{N}$, a set of branches $\mathcal{L}'$ such that ${\bf A}_{\mathcal{L}', \mathcal{N}'}={\bf A}_{\mathcal{L}, \mathcal{N}}$ and ${\bf Y}_{\mathcal{L}'}= {\bf B}_{\mathcal{L}}$, and a set of shunts $\mathcal{T}'$ such that ${\bf Y}_{\mathcal{T}'}= {\bf B}_{\mathcal{T}}$. The admittance matrix of the alternative network is:
\begin{equation}
{\bf Y}_{\mathcal{N}'} = {\bf A}_{\mathcal{L}, \mathcal{N}}^H {\bf B}_{\mathcal{L}} {\bf A}_{\mathcal{L}, \mathcal{N}} + {\bf B}_{\mathcal{T}},
\end{equation}
therefore
\begin{equation}
{\bf Y}_{\mathcal{N}} = jk {\bf Y}_{\mathcal{N}'}
\end{equation}
Notice that the alternative network satisfies Assumptions 1 and 2 and is purely resistive with no negative conductances. Hence ${\bf Y}_{\mathcal{N}'}$ satisfies Theorem 2. As $jk \neq 0$ we have that ${\rm Null} \left({{\bf Y}_{\mathcal{N}}}\right) = {\rm Null} \left({{\bf Y}_{\mathcal{N}'}}\right)$. Moreover, we know that ${\bf A}_{\mathcal{L}', \mathcal{N}'} = {\bf A}_{\mathcal{L}, \mathcal{N}}$, so ${\rm Null}\left({{\bf Y}_{\mathcal{N}} }\right) \subseteq {\rm Null} \left({{\bf A}_{\mathcal{L}, \mathcal{N}}}\right)$. $\hfill \square$

%We notice that if Condition \textit{2)} is satisfied then the network has no connection to ground, and so its equivalent admittance from any node to ground is zero. This means that Conditions \textit{1)} and \textit{2)} are mutually exclusive. 

\begin{figure}[t]
\tikzstyle{startstop} = [rectangle, rounded corners, minimum width=3cm, minimum height=1cm,text centered, draw=black, text width=3.5cm, fill=blue!30]
\tikzstyle{io} = [trapezium, trapezium left angle=70, trapezium right angle=110, minimum width=0.5cm, minimum height=0.65cm, text centered, draw=black, fill=green!30]
\tikzstyle{process} = [rectangle, minimum height=0.65cm, text centered, draw=black, text width=3.5cm, fill=orange!15]
\tikzstyle{decision} = [diamond, minimum width=0.5cm, text centered, draw=black, aspect=1.5, text width=1.5cm, fill=cyan!20]
\tikzstyle{line} = [thick,-,>=stealth]
\tikzstyle{arrow} = [thick,->,>=stealth]

\scalebox{0.725}{
\begin{tikzpicture}[x=1cm,y=1cm,node distance=1.5cm]
\node (pro1) [process] {Sort lines and reduce parallels};
\node (in1) [io, right of=pro1, xshift=2.0cm, inner xsep=-0.1cm] {Get $(\mathcal{N,L,T})$};
\node (start) [startstop, right of=in1, xshift=2.0cm] {Start};
\node (dec1) [decision, below of=pro1, yshift=-0.5cm, minimum width=3.0cm, aspect=1.0] {} (dec1) node[text centered, text width=2.0cm] {Thm. 2 assumptions hold?};
\node (pro2) [process, below of=dec1, yshift=-0.5cm] {Compute $(\mathcal{N',L',T'})$};
\node (pro3) [process, below of=pro2] {Get all reac. comps. $(\mathcal{N}'(k),\mathcal{L}'(k),\mathcal{T}'(k))$ for $k=1,\cdots,K$};
\node (pro4) [process, below of=pro3] {Set $k \gets 1$};
\node (dec2) [decision, below of=pro4, yshift=0.0cm, minimum width=2.0cm, aspect=1.5] {} (dec2) node[text centered] {$k \leq K$?};
\node (pro5) [process, below of=dec2, text width=1.6cm] {$k \gets k+1$};
\node (dec3) [decision, below of=pro5, yshift=-0.75cm, minimum width=4.5cm, aspect=0.8] {} (dec3) node[text centered, text width=3.0cm, yshift=-0.1cm] {Thm. 3 holds for $(\mathcal{N}'(k),\mathcal{L}'(k),\mathcal{T}'(k))$ ?};
\node (dec4) [decision, right of=dec2, xshift=2.0cm, minimum width=2.0cm, aspect=1.5] {} (dec4) node[text centered] {$\mathcal{T} = \emptyset$?};
\node (dec5) [above of=dec4, yshift=0.5cm] {};
\node (out1) [io, right of=dec5, xshift=2.0cm, inner xsep=-0.3cm] {``${\bf Y}_\mathcal{N}$ is invertible''};
\node (out2) [io, above of=dec5, yshift=0.5cm, text width=2.0cm, inner xsep=0.0cm] {``${\rm rank}\left({{\bf Y}_\mathcal{N}}\right) =$  $ {\rm rank}\left({{\bf A}_{\mathcal{L},\mathcal{N}}}\right)$''};
\draw (out2 -| out1) node(end) [startstop] {End};
\draw (dec1 -| end) node(out3) [io, text width=2.0cm, inner xsep=-0.15cm] {``Theorems cannot be applied''};

\draw [arrow] (start) -- (in1);
\draw [arrow] (in1) -- (pro1);
\draw [arrow] (pro1) -- (dec1);
\draw [arrow] (pro2) -- (pro3);
\draw [arrow] (pro3) -- (pro4);
\draw [arrow] (pro4) -- (dec2);
\draw [arrow] (out1) -- (end);
\draw [arrow] (out2) -- (end);
\draw [arrow] (out3) -- (end);
\draw [arrow] (pro5) -- (dec2);

\draw [arrow] (dec1) -- node[anchor=east] {yes} (pro2);
\draw [arrow] (dec1) -- node[anchor=south] {no} (out3);
\draw [arrow] (dec2) -- node[anchor=north] {no} (dec4);
\draw [arrow] (dec3) -- node[anchor=west] {yes} (pro5);
\draw [arrow] (dec4) -- node[anchor=east] {yes} (out2);
%\draw [arrow] (dec5) -- node[anchor=south] {yes} (out1);
%\draw [arrow] (dec5) -- node[anchor=east] {no} (out2);

\draw [arrow] (dec2) -- node[anchor=south] {yes} ++(-2.75,0.0) |- (dec3);
\draw [arrow] (dec4) -- node[anchor=north] {no} (dec4 -| out1) -- (out1);
\draw [arrow] (dec3) -- node[anchor=north] {no} ($(dec3 -| out3) + (2.25,0.0)$) |- (out3);
\end{tikzpicture}
}

\caption{\textcolor{\mycolor}{Flowchart describing an algorithm to certify the invertibility (or singularity) of an admittance matrix through the use of Theorems 1 to 3.}}
\label{fig:algorithm}
\vspace{-1em}
\end{figure}
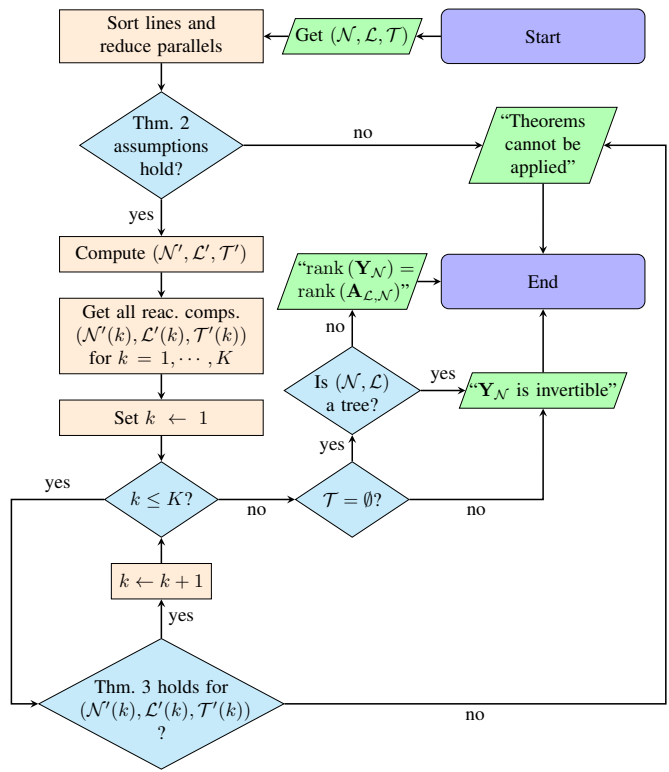

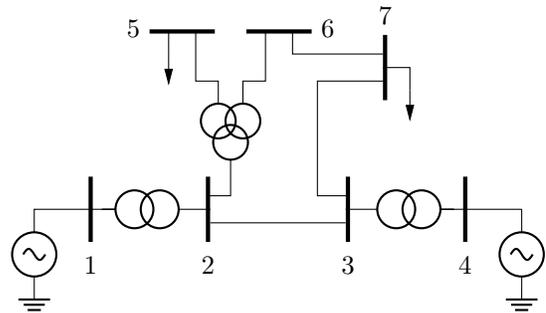
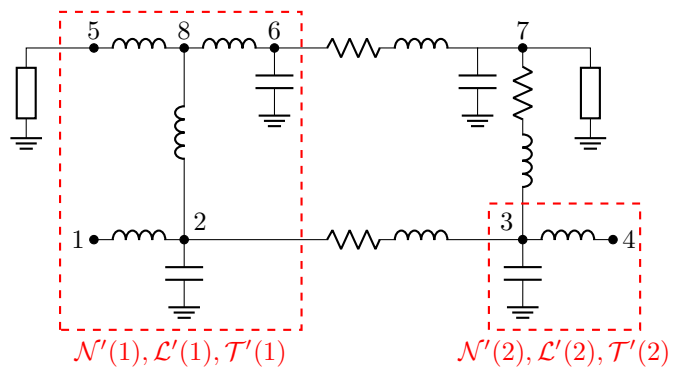
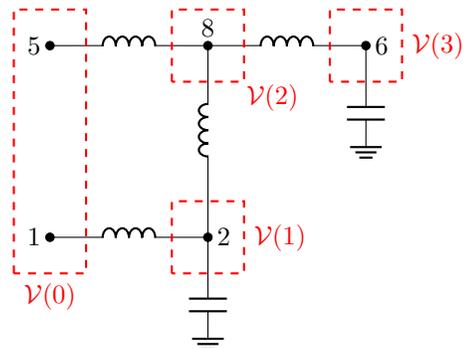
\begin{figure}[t]
\centering
\vspace{-1em}
\subfloat[One-line diagram of the example system.]{
\scalebox{1.0}{
\begin{circuitikz}[american voltages, scale=0.6]
 \ctikzset{monopoles/vee/arrow={Triangle[width=0.4*\scaledwidth, length=0.8*\scaledwidth]}}
\draw
  (0,0)
  node[tlground]{}
  to[sV, sources/scale=0.7, fill=white] ++(0,2) -- ++(0.5,0)
  to[bus=0.6, l_=$1$] ++(1.5,0)
  ++(-0.5,0) to[oosourcetrans, fill=white] ++(2,0)
  -- ++(0.3,0)
  to[bus=0.6, l_=$2$, name=n2] ++(0.1,0)
  ++(0,-0.3)-- ++(3,0) ++(0,0.3)
  to[bus=0.6, l_=$3$, name=n3] ++(0.1,0)
  -- ++(0.3,0)
  to[oosourcetrans, fill=white] ++(2,0)
  ++(-0.5,0) to[bus=0.6, l_=$4$] ++(1.5,0) -- ++(0.5,0)
  ++(0,-2) node[tlground]{}
  to[sV, sources/scale=0.7, fill=white] ++(0,2)
  ;
\draw[color=white]
  (n2.center) ++(0.5,0.3) to[ooosource, name=twt, color=black, fill=white] ++(0,3)
  ;
\draw
  (n2.center) ++(0,0.3) -- ++(0.5,0) -- (twt.left)
  (twt.sec1) -- ++(0,0.5) -- ++(-0.5,0) -- ++(0,1.1)
  ++(-0.3,0) to[bus=0.6, l=$5$] ++(0,0)
  ++(-0.3,0) -- ++(0,-0.5) node[vee]{}
  (twt.tert1) -- ++(0,0.5) -- ++(0.5,0) -- ++(0,1.1)
  ++(0.3,0) to[bus=0.6, l_=$6$] ++(0,0)
  ++(0.3,0) -- ++(0,-0.5) -- ++(2,0)
  ++(0,-0.3) to[bus=0.6, l=$7$, name=n7] ++(0.1,0)
  -- ++(0.5,0) -- ++(0,-0.5) node[vee]{}
  (n7.center) ++(0,-0.3) -- ++(-1.5,0) |- ($(n3.center) + (0,0.3)$)
  ;
\end{circuitikz}
} \label{fig:alg_example_circuit}
} \\
\subfloat[Circuit model of the system with its reactive components outlined.]{
\scalebox{1.0}{
\begin{circuitikz}[american voltages, american resistors, scale=0.6]
\ctikzset{resistors/scale=0.6, inductors/scale=0.6, capacitors/scale=0.6}
\draw
  (0,0) node[name=n1]{} node[anchor=east]{$1$}
  to[L, *-*] ++(2,0) node[name=n2]{} node[anchor=south west]{$2$}
  to[C] ++(0,-1.5) node[tlground]{}
  (n2.center) -- ++(0,0.5) to[L, -*] ++(0,3.75) node[name=n8]{} node[anchor=south]{$8$}
  to[L, -*, mirror] ++(-2,0) node[anchor=south]{$5$}
  -- ++(-1.5,0) to[R, european resistors] ++(0,-2) node[tlground]{}
  (n8.center) to[L, -*] ++(2,0) node[name=n6]{} node[anchor=south]{$6$}
  to[C] ++(0,-1.5) node[tlground]{}
  (n6.center) -- ++(1,0) to[R] ++(1.5,0)
  to[L] ++(1.5,0) -- ++(0.5,0) 
  to[C] ++(0,-1.5) node[tlground]{} ++(0,1.5)
  -- ++(1,0) node[circ, name=n7]{} node[anchor=south]{$7$}
  -- ++(1.5,0) to[R, european resistors] ++(0,-2) node[tlground]{}
  (n7.center) -- ++(0,-0.25) to[R] ++(0,-1.5)
  to[L] ++(0,-1.5) -- ++(0,-1) node[circ, name=n3]{} node[anchor=south east]{$3$}
  (n2.center) -- ++(3,0) to[R] ++(1.5,0)
  to[L] ++(1.5,0) -- (n3.center)
  to[C] ++(0,-1.5) node[tlground]{}
  (n3.center) to[L] ++(2,0) node[circ]{} node[name=n4]{} node[anchor=west]{$4$}
  ;
  \draw[red, thick, dashed] 
  ($(n1.center)+(-0.75,-2)$) node[name=rec1_sw]{}
  rectangle ($(n6.center)+(0.6,0.8)$) node[midway, name=rec1_center]{}
  (rec1_center |- rec1_sw) node[anchor=north]{$\mathcal{N}'(1), \mathcal{L}'(1), \mathcal{T}'(1)$}
  ($(n3.center)+(-0.75,-2)$) node[name=rec2_sw]{}
  rectangle ($(n4.center)+(0.6,0.8)$) node[midway, name=rec2_center]{}
  (rec2_center |- rec2_sw) node[anchor=north]{$\mathcal{N}'(2), \mathcal{L}'(2), \mathcal{T}'(2)$}
  ;
\end{circuitikz}
} \label{fig:alg_example_comps}
} \\
\subfloat[Node partition for the tree $(\mathcal{N}'(1),\mathcal{L}'(1))$, using node 6 as root.]{
\scalebox{1.0}{
\begin{circuitikz}[american voltages, american resistors, scale=0.6]
\ctikzset{resistors/scale=0.6, inductors/scale=0.6, capacitors/scale=0.6}
\draw
  (0,0) node[circ, name=n1]{} node[anchor=east]{$1$}
  to[L] ++(3.5,0)
  node[circ, name=n2]{} node[anchor=west]{$2$}
  -- ++(0,-0.75) to[C] ++(0,-1.5) node[tlground]{}
  (n2.center) -- ++(0,0.5) to[L, -*] ++(0,3.75) node[name=n8]{} node[anchor=south]{$8$}
  to[L, mirror] ++(-3.5,0)
  node[circ, name=n5]{} node[anchor=east]{$5$}
  (n8.center) to[L, -*] ++(3.5,0) node[name=n6]{} node[anchor=west]{$6$}
  -- ++(0,-0.75) to[C] ++(0,-1.5) node[tlground]{}
  ;
  \draw[red, thick, dashed]
  ($(n1.center)+(-0.8,-0.8)$) node[name=rec0_sw]{}
  rectangle ($(n5.center)+(0.8,0.8)$) node[midway, name=rec0_center]{}
  (rec0_center |- rec0_sw) node[anchor=north]{$\mathcal{V}(0)$}
  ($(n2.center)+(-0.8,-0.8)$) rectangle ($(n2.center)+(0.8,0.8)$) node[midway, name=rec1_center]{} node[name=rec1_ne]{}
  (rec1_center -| rec1_ne) node[anchor=west]{$\mathcal{V}(1)$}
  ($(n8.center)+(-0.8,-0.8)$) node[name=rec2_sw]{}
  rectangle ($(n8.center)+(0.8,0.8)$) node[name=rec2_ne]{}
  (rec2_sw -| rec2_ne) node[anchor=north west, shift=({-0.1,0.1})]{$\mathcal{V}(2)$}
  ($(n6.center)+(-0.8,-0.8)$) rectangle ($(n6.center)+(0.8,0.8)$) node[midway, name=rec3_center]{} node[name=rec3_ne]{}
  (rec3_center -| rec3_ne) node[anchor=west]{$\mathcal{V}(3)$}
  ;
\end{circuitikz}
} \label{fig:alg_example_tree}
}%
\caption{Example system to illustrate how to apply the main theorems.}
\label{fig:alg_example}
\vspace{-1em}
\end{figure}

\textcolor{\mycolor}{We now have enough tools to construct an algorithm to check the invertibility of the admittance matrix. First, we reduce any parallel lines in order to comply with Remark~1. Then we check if the network satisfies the assumptions of Theorem~2; if so, we compute the reactive subsystem $(\mathcal{N',L',T'})$ by removing all elements with positive resistance. Afterwards, we compute all the $K$ connected components of $(\mathcal{N',L',T'})$ (when $(\mathcal{N',L',T'})$ is empty, then $K=0$). These connected components are computed using the Breadth First Search (BFS) algorithm \cite{algorithms}, whose complexity is $\mathcal{O}\left({|\mathcal{N}| + |\mathcal{L}|}\right)$ (linear in the system size). For each connected component $(\mathcal{N}'(k),\mathcal{L}'(k),\mathcal{T}'(k)), k=1,\ldots,K$, we check if Theorem~3 can be applied to the component through any of its conditions. If Theorem~3 holds for all reactive components, then Theorem~2 holds for the networks and thus Theorem~1 holds as well. Finally, using Theorem~1, we can certify the invertibility of the admittance matrix if the network has shunts. Otherwise, ${\rm rank}\left({{\bf A}_{\mathcal{L},\mathcal{N}}}\right)$ needs to be computed. In the special case that the network is radial, we have that ${\rm rank}\left({{\bf A}_{\mathcal{L},\mathcal{N}}}\right) = |\mathcal{N}| - 1$, and thus if there are no shunts we can certify that the admittance matrix is \textit{singular}. A flowchart of the algorithm is shown Fig.~\ref{fig:algorithm}. To illustrate the idea behind the algorithm, consider the example system of Fig.~\ref{fig:alg_example}. The one-line diagram of the system is shown in Fig.~\ref{fig:alg_example_circuit}. In Fig.~\ref{fig:alg_example_comps}, we have the circuit model of the system, where the loads are modelled as constant admittances and each transmission line is modeled using a $\pi$ circuit. The example system possesses two reactive components, outlined in the figure (shunt loads are not included in the components, as they have a resistive part). If the main condition of Theorem~2 can be proved for each component (by means of Theorem~3, for example), then Theorem~2 will hold for the system, and thus Theorem~1 will holds for the system as well. The branches of the first component, $(\mathcal{N}'(1),\mathcal{T}'(1))$, form a tree. Choosing node 6 as root, we obtain the node partition shown in Fig.~\ref{fig:alg_example_tree}. This partition can be used to check Condition \textit{1)} of Theorem~3.}

\section{Implementation and Test Cases}
\label{sec:implementation}

\begin{table}[t]
\centering
\caption{PGLib Test Cases Used for Checking the Theorems}
\label{tbl:test_cases}
\hspace*{-1em}
\textcolor{\mycolor}{
\scalebox{0.95}{
\begin{tabular}{|l|r|r|r|l|l|}
\hline
Test case & $|\mathcal{N}|$ & $|\mathcal{L}|$ & Reac. & Satisfy thm. & ${\bf Y}_{\mathcal{N}}$ non- \\
 & & & line \% & conditions? & singular? \\
\hline
case3\_lmbd            &3      &3      & 0.0\%    &Yes    &Yes   \\
\hline
case5\_pjm             &5      &6      & 0.0\%    &Yes    &Yes   \\
\hline
case14\_ieee           &14     &20     &25.0\%    &No     &-     \\
\hline
case24\_ieee\_rts       &24     &38     & 0.0\%    &Yes    &Yes   \\
\hline
case30\_as             &30     &41     &17.1\%    &No     &-     \\
\hline
case30\_ieee           &30     &41     &17.1\%    &No     &-     \\
\hline
case39\_epri           &39     &46     & 8.7\%    &Yes    &Yes   \\
\hline
case57\_ieee           &57     &80     &22.5\%    &No     &-     \\
\hline
case60\_c              &60     &88     &40.9\%    &Yes    &Yes   \\
\hline
case73\_ieee\_rts       &73     &120    & 0.8\%    &Yes    &Yes   \\
\hline
case89\_pegase         &89     &210    & 4.8\%    &Yes    &Yes   \\
\hline
case118\_ieee          &118    &186    & 4.8\%    &Yes    &Yes   \\
\hline
case162\_ieee\_dtc      &162    &284    &11.6\%    &Yes    &Yes   \\
\hline
case179\_goc           &179    &263    &27.4\%    &Yes    &Yes   \\
\hline
case200\_activ         &200    &245    & 0.0\%    &Yes    &Yes   \\
\hline
case240\_pserc         &240    &448    &20.8\%    &Yes    &Yes   \\
\hline
case300\_ieee          &300    &411    &15.6\%    &No     &-     \\
\hline
case500\_goc           &500    &733    & 0.0\%    &Yes    &Yes   \\
\hline
case588\_sdet          &588    &686    & 7.1\%    &No     &-     \\
\hline
case793\_goc           &793    &913    & 0.7\%    &No     &-     \\
\hline
case1354\_pegase       &1354   &1991   & 0.1\%    &Yes    &Yes   \\
\hline
case1888\_rte          &1888   &2531   & 9.8\%    &Yes    &Yes   \\
\hline
case1951\_rte          &1951   &2596   &13.1\%    &Yes    &Yes   \\
\hline
case2000\_goc          &2000   &3639   & 0.0\%    &Yes    &Yes   \\
\hline
case2312\_goc          &2312   &3013   & 0.0\%    &No     &-     \\
\hline
case2383wp\_k          &2383   &2896   & 6.7\%    &Yes    &Yes   \\
\hline
case2736sp\_k          &2736   &3504   & 1.1\%    &Yes    &Yes   \\
\hline
case2737sop\_k         &2737   &3506   & 1.1\%    &No     &-     \\
\hline
case2742\_goc          &2742   &4673   & 0.0\%    &Yes    &Yes   \\
\hline
case2746wop\_k         &2746   &3514   & 1.1\%    &No     &-     \\
\hline
case2746wp\_k          &2746   &3514   & 1.1\%    &Yes    &Yes   \\
\hline
case2848\_rte          &2848   &3776   & 5.5\%    &Yes    &Yes   \\
\hline
case2853\_sdet         &2853   &3921   & 9.7\%    &No     &-     \\
\hline
case2868\_rte          &2868   &3808   & 6.7\%    &Yes    &Yes   \\
\hline
case2869\_pegase       &2869   &4582   & 3.0\%    &No     &-     \\
\hline
case3012wp\_k          &3012   &3572   & 0.3\%    &No     &-     \\
\hline
case3022\_goc          &3022   &4135   & 5.2\%    &No     &-     \\
\hline
case3120sp\_k          &3120   &3693   & 0.3\%    &No     &-     \\
\hline
case3375wp\_k          &3374   &4161   & 0.6\%    &No     &-     \\
\hline
case3970\_goc          &3970   &6641   & 0.0\%    &Yes    &Yes   \\
\hline
case4020\_goc          &4020   &6988   & 0.0\%    &Yes    &Yes   \\
\hline
case4601\_goc          &4601   &7199   & 0.0\%    &Yes    &Yes   \\
\hline
case4619\_goc          &4619   &8150   & 0.0\%    &Yes    &Yes   \\
\hline
case4661\_sdet         &4661   &5997   & 1.6\%    &No     &-     \\
\hline
case4837\_goc          &4837   &7765   & 0.0\%    &Yes    &Yes   \\
\hline
case4917\_goc          &4917   &6726   & 2.6\%    &No     &-     \\
\hline
case6468\_rte          &6468   &9000   & 2.1\%    &Yes    &Yes   \\
\hline
case6470\_rte          &6470   &9005   & 2.4\%    &Yes    &Yes   \\
\hline
case6495\_rte          &6495   &9019   & 2.8\%    &Yes    &Yes   \\
\hline
case6515\_rte          &6515   &9037   & 2.9\%    &Yes    &Yes   \\
\hline
case8387\_pegase       &8387   &14561  & 4.3\%    &No     &-     \\
\hline
case9241\_pegase       &9241   &16049  & 5.3\%    &No     &-     \\
\hline
case9591\_goc          &9591   &15915  & 0.0\%    &Yes    &Yes   \\
\hline
case10000\_goc         &10000  &13193  & 0.0\%    &Yes    &Yes   \\
\hline
case10480\_goc         &10480  &18559  & 0.0\%    &Yes    &Yes   \\
\hline
case13659\_pegase      &13659  &20467  & 7.0\%    &No     &-     \\
\hline
case19402\_goc         &19402  &34704  & 0.0\%    &Yes    &Yes   \\
\hline
case24464\_goc         &24464  &37816  & 0.0\%    &Yes    &Yes   \\
\hline
case30000\_goc         &30000  &35393  & 0.0\%    &Yes    &Yes   \\
\hline
\end{tabular}
}
}
\vspace*{-1em}
\end{table}

We developed MATLAB R2012b code that implements the algorithm described in Section~\ref{sec:main_results}. The code is publicly available at the following page:

\phantom{.}

\noindent
\href{https://github.com/djturizo/ybus-inv-check}{\textcolor{blue}{\texttt{https://github.com/djturizo/ybus-inv-check}}}

\phantom{.}

\noindent
This code is not optimized for performance, but rather serves as a proof-of-concept for the complexity of the algorithm. The interested reader can examine the code and its comments to see that the program has a time complexity of $\mathcal{O}\left({|\mathcal{N}| + |\mathcal{L}|}\right)$ (linear in the system size)\footnote{Our implementation relies on standard algorithms like counting sort, radix sort and BFS. These algorithms are known to have linear time complexity \cite{algorithms}, so the whole implementation has linear complexity as well.}. We remark that comparisons cannot be exact due to the finite-precision computations, so the program uses a user-defined tolerance for all comparisons. 
%Directly computing the rank of the generalized incidence matrix ${\bf A}_{\mathcal{L},\mathcal{N}}$ cannot be done in linear time, so the program does not compute the rank of the admittance matrix when there are no shunts.

%For the numerical experiments, we selected the test cases from the Power Grid Library PGLib \cite{pglib} (from the OPF benchmarks, more specifically). Some of the PGLib test cases have a small number of negative resistance elements, precluding the application of Theorem~1. This is the result of modeling choices associated with equivalenced networks~\cite{josz_pegase}. Since these non-passive branches are of an artificial (non-physical) nature, we focus on the other 44 PGLib test cases without negative resistance elements for our numerical experiments.

\textcolor{\mycolor}{For the numerical experiments we ran the program using radial (distribution) and meshed (transmission) test cases. The radial test were taken from M{\sc atpower} \cite{matpower_manual}. The program successfully applied the theorems to certify the invertibility (or singularity) of the admittance matrix for all the readial test cases of M{\sc atpower}. The meshed test cases were selected from the Power Grid Library PGLib \cite{pglib} (from the OPF benchmarks, more specifically).} Some of the PGLib test cases have a small number of negative resistance elements, precluding the application of Theorem~1. This is the result of modeling choices associated with equivalenced networks~\cite{josz_pegase}. Since these non-passive branches are of an artificial (non-physical) nature, we focus on the other 44 PGLib test cases without negative resistance elements for our numerical experiments. With a tolerance of $10^{-12}$, we obtained the results shown in Table~\ref{tbl:test_cases}. \textcolor{\mycolor}{Note that the fourth column of this table refers to the percentage of branches in the system that are purely reactive.} Of the 44 test cases, we found that 6 of them did not satisfy the conditions of the theorems. Thus, the program could not certify the invertibility of the admittance matrix for those 6 cases. These cases are identified with a dash in the last column of Table~\ref{tbl:test_cases} to indicate that the theorems cannot certify whether or not the admittance matrix is invertible. The reason why the invertibility could not be certified for each of the 6 cases is because they have reactive components with topologies not covered by Theorem~3. Such components have inductors, capacitors and loops formed by branches. However, those complex topologies are uncommon, as for the other 38 cases (86\% of all cases) the conditions of the theorems hold, so the program can certify the whether the admittance matrix is invertible or not for each case. The admittance matrix is known to be invertible for the test cases, so we get positive results whenever the theorems were applicable. \textcolor{\mycolor}{In contrast, none of the 41 cases satisfy the conditions of the invertibility theorems developed in \cite{rank_1phase, rank_3phase, gatsis, low_theorem}. That is, while none of the existing theoretical results regarding the invertibility of the admittance matrix can be applied to any of the considered PGLib test cases, our results successfully certify the invertibility of the admittance matrix in 86\% of these test cases.} The results show that the theorems can be used to certify the invertibility of the admittance matrices for a wide range of practical and realistic power systems. Moreover, the cases where the theorems cannot be applied to a realistic power system are uncommon.

\section{Conclusions}

% Summarize what we did, then discuss future work.

This paper studied the invertibility of the admittance matrix for balanced networks. First, we analyzed a theorem from the literature regarding conditions guaranteeing invertibility of the admittance matrix, and we found a technical issue in the proof of that theorem. Next, we developed a framework of lemmas and assumptions that allowed us to amend the proof of previous claims, developing relaxed conditions that guarantee the invertibility of the admittance matrix and generalizing the results to systems with branches modeled as purely reactive elements and transformers with off-nominal tap ratios. Finally, we implemented and publicly released a proof-of-concept program that uses the theorems to certify the invertibility of the admittance matrix. Numerical tests showed that the theorems are applicable in a large number of realistic power systems.
%This paper studied the invertibility of the admittance matrix for balanced networks. First, we analyzed claims from the literature regarding conditions guaranteeing invertibility of the admittance matrix. We found a flaw in these claims and provided a counterexample. Next, we developed a framework of lemmas and assumptions that allowed us to amend the flaws of previous claims, developing various conditions that guarantee the invertibility of the admittance matrix and generalizing the results to systems with phase-shifting transformers. Finally, we discussed these conditions' implications, their wide practical applicability, and their usefulness in enforcing invertibility.

The theory developed in this paper has solely considered admittance matrices for balanced single-phase equivalent network representations. With rapidly increasing penetration of distributed energy resources in unbalanced distribution systems, extending the theory developed here to address the admittance matrices associated with polyphase networks is an important direction for future work. The authors of~\cite{rank_1phase} considered this topic in~\cite{rank_3phase}, where they generalize Theorem~1 to polyphase networks. However, the theory in~\cite{rank_3phase} also relies on the incorrectly stated Lemma~3 and hence may also benefit from amendments and extensions similar to those in this paper.
%We believe that similar amendments to those in this paper could correct and extend the results in~\cite{rank_3phase}.

\section*{Acknowledgements}
The authors greatly appreciate technical discussions with Mario Paolone and Andreas Kettner.

% if have a single appendix:
%\appendix[Proof of the Zonklar Equations]
% or
%\appendix  % for no appendix heading
% do not use \section anymore after \appendix, only \section*
% is possibly needed

% use appendices with more than one appendix
% then use \section to start each appendix
% you must declare a \section before using any
% \subsection or using \label (\appendices by itself
% starts a section numbered zero.)
%

% Can use something like this to put references on a page
% by themselves when using endfloat and the captionsoff option.
\ifCLASSOPTIONcaptionsoff
  \newpage
\fi

% trigger a \newpage just before the given reference
% number - used to balance the columns on the last page
% adjust value as needed - may need to be readjusted if
% the document is modified later
%\IEEEtriggeratref{8}
% The "triggered" command can be changed if desired:
%\IEEEtriggercmd{\enlargethispage{-5in}}

% references section

% can use a bibliography generated by BibTeX as a .bbl file
% BibTeX documentation can be easily obtained at:
% http://mirror.ctan.org/biblio/bibtex/contrib/doc/
% The IEEEtran BibTeX style support page is at:
% http://www.michaelshell.org/tex/ieeetran/bibtex/
\bibliographystyle{IEEEtran}
% argument is your BibTeX string definitions and bibliography database(s)
\bibliography{IEEEabrv,refs}

% Generated by IEEEtran.bst, version: 1.14 (2015/08/26)
\begin{thebibliography}{10}
\providecommand{\url}[1]{#1}
\csname url@samestyle\endcsname
\providecommand{\newblock}{\relax}
\providecommand{\bibinfo}[2]{#2}
\providecommand{\BIBentrySTDinterwordspacing}{\spaceskip=0pt\relax}
\providecommand{\BIBentryALTinterwordstretchfactor}{4}
\providecommand{\BIBentryALTinterwordspacing}{\spaceskip=\fontdimen2\font plus
\BIBentryALTinterwordstretchfactor\fontdimen3\font minus
  \fontdimen4\font\relax}
\providecommand{\BIBforeignlanguage}[2]{{%
\expandafter\ifx\csname l@#1\endcsname\relax
\typeout{** WARNING: IEEEtran.bst: No hyphenation pattern has been}%
\typeout{** loaded for the language `#1'. Using the pattern for}%
\typeout{** the default language instead.}%
\else
\language=\csname l@#1\endcsname
\fi
#2}}
\providecommand{\BIBdecl}{\relax}
\BIBdecl

\bibitem{kundur}
P.~Kundur, \emph{Power System Stability And Control}.\hskip 1em plus 0.5em
  minus 0.4em\relax McGraw-Hill, 1994.

\bibitem{crow}
M.~L. Crow, \emph{Computational Methods for Electric Power Systems},
  3rd~ed.\hskip 1em plus 0.5em minus 0.4em\relax CRC Press, 2015.

\bibitem{kron}
G.~Kron, \emph{Tensors for Circuits}.\hskip 1em plus 0.5em minus 0.4em\relax
  Dover Publications, 1959.

\bibitem{efficient_thevenin}
C.~H.~L. {Jørgensen}, J.~G. {Møller}, S.~{Sommer}, and H.~{Jóhannsson}, ``{A
  Memory-Efficient Parallelizable Method for Computation of Thévenin
  Equivalents Used in Real-Time Stability Assessment},'' \emph{IEEE
  Transactions on Power Systems}, vol.~34, no.~4, pp. 2675--2684, 2019.

\bibitem{impedance_estimator}
P.~{Cuffe} and F.~{Milano}, ``{Validating Two Novel Equivalent Impedance
  Estimators},'' \emph{IEEE Transactions on Power Systems}, vol.~33, no.~1, pp.
  1151--1152, 2018.

\bibitem{frequency_divider}
F.~{Milano} and {\'{A}}.~{Ortega}, ``{Frequency Divider},'' \emph{IEEE
  Transactions on Power Systems}, vol.~32, no.~2, pp. 1493--1501, 2017.

\bibitem{anderson_faults}
P.~M. Anderson, \emph{Analysis of Faulted Power Systems}.\hskip 1em plus 0.5em
  minus 0.4em\relax Wiley-IEEE, 1995.

\bibitem{stott2009}
B.~Stott, J.~Jardim, and O.~Alsa\c{c}, ``{DC Power Flow Revisited},''
  \emph{IEEE Transactions on Power Systems}, vol.~24, no.~3, pp. 1290--1300,
  2009.

\bibitem{ptdf}
M.~K. {Enns}, J.~J. {Quada}, and B.~{Sackett}, ``{Fast Linear Contingency
  Analysis},'' \emph{IEEE Transactions on Power Apparatus and Systems}, vol.
  PAS-101, no.~4, pp. 783--791, 1982.

\bibitem{atc}
G.~C. {Ejebe}, J.~G. {Waight}, M.~{Sanots-Nieto}, and W.~F. {Tinney}, ``{Fast
  Calculation of Linear Available Transfer Capability},'' \emph{IEEE
  Transactions on Power Systems}, vol.~15, no.~3, pp. 1112--1116, 2000.

\bibitem{tenti2018}
P.~Tenti and T.~Caldognetto, ``{Optimal control of Local Area Energy Networks
  {(E-LAN)}},'' \emph{Sustainable Energy, Grids and Networks}, vol.~14, pp.
  12--24, 2018.

\bibitem{rahman2019}
M.~Rahman, V.~Cecchi, and K.~Miu, ``{Power Handling Capabilities of
  Transmission Systems using a Temperature-Dependent Power Flow},''
  \emph{Electric Power Systems Research}, vol. 169, pp. 241--249, 2019.

\bibitem{rrlu_fact}
L.~Miranian and M.~Gu, ``{Strong Rank Revealing {LU} Factorizations},''
  \emph{Linear Algebra and Its Applications}, vol. 367, pp. 1--16, 2003.

\bibitem{rrqr_fact}
M.~Gu and S.~C. Eisenstat, ``{Efficient Algorithms for Computing a Strong
  Rank-Revealing {QR} Factorization},'' \emph{SIAM Journal on Scientific
  Computing}, vol.~17, no.~4, pp. 848--869, 1996.

\bibitem{cond}
C.~P. Br{\'a}s, W.~W. Hager, and J.~J. J{\'u}dice, ``{An Investigation of
  Feasible Descent Algorithms for Estimating the Condition Number of a
  Matrix},'' \emph{TOP}, vol.~20, no.~3, pp. 791--809, Oct. 2012.

\bibitem{cond_reliable}
\BIBentryALTinterwordspacing
H.~Avron, A.~Druinsky, and S.~Toledo, ``Spectral condition-number estimation of
  large sparse matrices,'' \emph{Numerical Linear Algebra with Applications},
  vol.~26, no.~3, 2019. [Online]. Available:
  \url{https://onlinelibrary.wiley.com/doi/abs/10.1002/nla.2235}
\BIBentrySTDinterwordspacing

\bibitem{transmission_switching}
E.~B. {Fisher}, R.~P. {O'Neill}, and M.~C. {Ferris}, ``{Optimal Transmission
  Switching},'' \emph{IEEE Transactions on Power Systems}, vol.~23, no.~3,
  2008.

\bibitem{reconfiguration_1}
T.~Ding, Y.~Lin, Z.~Bie, and C.~Chen, ``{A Resilient Microgrid Formation
  Strategy for Load Restoration Considering Master-Slave Distributed Generators
  and Topology Reconfiguration},'' \emph{Applied Energy}, vol. 199, pp.
  205--216, 2017.

\bibitem{reconfiguration_2}
Z.~{Guo}, Z.~{Zhou}, and Y.~{Zhou}, ``{Impacts of Integrating Topology
  Reconfiguration and Vehicle-to-Grid Technologies on Distribution System
  Operation},'' \emph{IEEE Transactions on Sustainable Energy}, vol.~11, no.~2,
  pp. 1023--1032, 2020.

\bibitem{rank_1phase}
A.~M. {Kettner} and M.~{Paolone}, ``{On the Properties of the Power Systems
  Nodal Admittance Matrix},'' \emph{IEEE Transactions on Power Systems},
  vol.~33, no.~1, pp. 1130--1131, 2018.

\bibitem{rank_3phase}
------, ``{On the Properties of the Compound Nodal Admittance Matrix of
  Polyphase Power Systems},'' \emph{IEEE Transactions on Power Systems},
  vol.~34, no.~1, pp. 444--453, 2019.

\bibitem{gatsis}
M.~{Bazrafshan} and N.~{Gatsis}, ``{Comprehensive Modeling of Three-Phase
  Distribution Systems via the Bus Admittance Matrix},'' \emph{IEEE
  Transactions on Power Systems}, vol.~33, no.~2, pp. 2015--2029, 2018.

\bibitem{low_theorem}
\BIBentryALTinterwordspacing
S.~H. Low. (2022, June) Power {S}ystems {A}nalysis: A {M}athematical
  {A}pproach. {L}ecture {N}otes for {E}{E}/{C}{S}/{E}{S}{T} 135. California
  Institute of Technology. [Online]. Available:
  \url{http://netlab.caltech.edu/book/book.html}
\BIBentrySTDinterwordspacing

\bibitem{bakshi}
U.~Bakshi and A.~Bakshi, \emph{Electrical Networks}.\hskip 1em plus 0.5em minus
  0.4em\relax Technical Pub., 2008.

\bibitem{pglib}
{IEEE PES Task Force on Benchmarks for Validation of Emerging Power System
  Algorithms}, ``{The Power Grid Library for Benchmarking {AC} Optimal Power
  Flow Algorithms},'' \emph{arXiv:1908.02788v2}, Jan. 2021.

\bibitem{stevenson}
W.~Stevenson and J.~Grainger, \emph{Power System Analysis}.\hskip 1em plus
  0.5em minus 0.4em\relax McGraw-Hill, 1994.

\bibitem{algebraic_graph}
C.~Godsil and G.~F. Royle, \emph{Algebraic Graph Theory}.\hskip 1em plus 0.5em
  minus 0.4em\relax Springer Science \& Business Media, 2001, vol. 207.

\bibitem{arrillaga}
J.~Arrillaga and C.~Arnold, \emph{Computer Analysis of Power Systems}.\hskip
  1em plus 0.5em minus 0.4em\relax John Wiley and Sons Inc., 1990.

\bibitem{linear_algebra}
C.~D. Meyer, \emph{Matrix Analysis and Applied Linear Algebra}.\hskip 1em plus
  0.5em minus 0.4em\relax SIAM, 2000, vol.~71.

\bibitem{graph_theory}
R.~Diestel, \emph{Graph Theory}, 5th~ed.\hskip 1em plus 0.5em minus 0.4em\relax
  Springer, 2017.

\bibitem{quadratic_forms}
I.~R. Shafarevich and A.~O. Remizov, \emph{Linear Algebra and Geometry}.\hskip
  1em plus 0.5em minus 0.4em\relax Springer Science \& Business Media, 2013.

\bibitem{algorithms}
T.~H. Cormen, C.~E. Leiserson, R.~L. Rivest, and C.~Stein, \emph{Introduction
  to Algorithms}, 4th~ed.\hskip 1em plus 0.5em minus 0.4em\relax MIT Press,
  2022.

\bibitem{matpower_manual}
\BIBentryALTinterwordspacing
R.~D. Zimmerman and C.~E. Murillo-Sánchez, ``Matpower user's manual,'' 2020.
  [Online]. Available: \url{https://matpower.org/docs/MATPOWER-manual-7.1.pdf}
\BIBentrySTDinterwordspacing

\bibitem{josz_pegase}
C.~Josz, S.~Fliscounakis, J.~Maeght, and P.~Panciatici, ``{AC Power Flow Data
  in MATPOWER and QCQP Format: iTesla, RTE Snapshots, and PEGASE},''
  \emph{arXiv:1603.01533v3}, Mar. 2016.

\end{thebibliography}
%
% <OR> manually copy in the resultant .bbl file
% set second argument of \begin to the number of references
% (used to reserve space for the reference number labels box)

% biography section
% 
% If you have an EPS/PDF photo (graphicx package needed) extra braces are
% needed around the contents of the optional argument to biography to prevent
% the LaTeX parser from getting confused when it sees the complicated
% \includegraphics command within an optional argument. (You could create
% your own custom macro containing the \includegraphics command to make things
% simpler here.)
%\begin{IEEEbiography}[{\includegraphics[width=1in,height=1.25in,clip,keepaspectratio]{mshell}}]{Michael Shell}
% or if you just want to reserve a space for a photo:

% \begin{IEEEbiography}{Michael Shell}
% Biography text here.
% \end{IEEEbiography}

\begin{IEEEbiography}[{\includegraphics[width=1in,height=1.25in,clip,keepaspectratio]{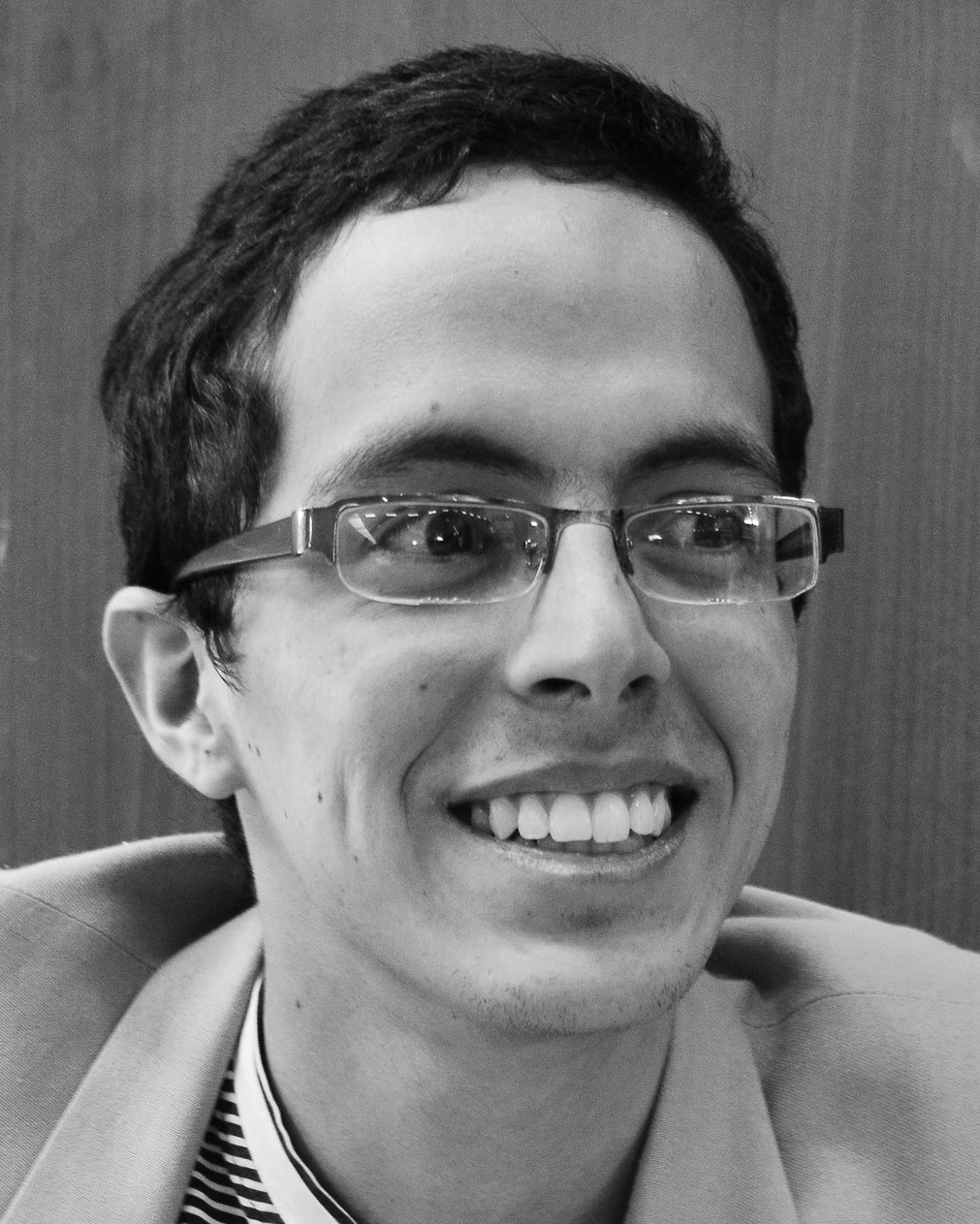}}]{Daniel Turizo}
(M'20) received the B.S. and M.S. degrees in Electrical Engineering from the Universidad del Norte, Barranquilla, Colombia in 2016 and 2018, respectively. He is currently a Ph.D. student at the School of Electrical and Computer Engineering, Georgia Institute of Technology, Atlanta, GA, USA where he is a Fullbright Fellow. From 2018-2020 he was an Adjunct Professor of Electrical Engineering at the Universidad del Norte, Barranquilla, Colombia. 
\end{IEEEbiography}

\begin{IEEEbiography}[{\includegraphics[width=1in,height=1.25in,clip,keepaspectratio]{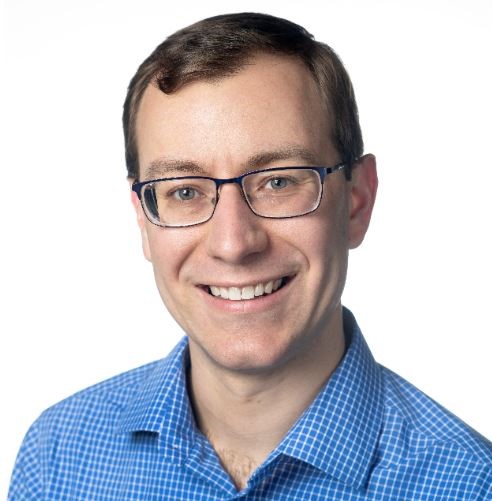}}]{Daniel K. Molzahn}
(S’09-M’13-SM’19) is an Assistant Professor in the School of Electrical and Computer Engineering at the Georgia Institute of Technology and also holds an appointment as a computational engineer in the Energy Systems Division at Argonne National Laboratory. He was a Dow Postdoctoral Fellow in Sustainability at the University of Michigan, Ann Arbor. He received the B.S., M.S., and Ph.D. degrees in electrical engineering and the Masters of Public Affairs degree from the University of Wisconsin--Madison, where he was a National Science Foundation Graduate Research Fellow. He received the IEEE Power and Energy Society's Outstanding Young Engineer Award in 2021 and the National Science Foundation's CAREER Award in 2022.
\end{IEEEbiography}

% % if you will not have a photo at all:
% \begin{IEEEbiographynophoto}{John Doe}
% Biography text here.
% \end{IEEEbiographynophoto}

% % insert where needed to balance the two columns on the last page with
% % biographies
% %\newpage

% \begin{IEEEbiographynophoto}{Jane Doe}
% Biography text here.
% \end{IEEEbiographynophoto}

% You can push biographies down or up by placing
% a \vfill before or after them. The appropriate
% use of \vfill depends on what kind of text is
% on the last page and whether or not the columns
% are being equalized.

%\vfill

% Can be used to pull up biographies so that the bottom of the last one
% is flush with the other column.
%\enlargethispage{-5in}

% that's all folks
\end{document}